%% file: Main.tex
\documentclass[a4paper,9pt]{article}
\usepackage[left=1.9cm,right=1.9cm,top=2.54cm,bottom=2.54cm]{geometry} %Advances in Mathematics
\usepackage{indentfirst}
\usepackage{amscd}
\usepackage{amsmath}
\usepackage{amssymb}
\usepackage{amsthm} 

\usepackage{pinlabel} %Pin­la­bel is a la­belling pack­age for at­tach­ing per­fectly for­mat­ted TeX la­bels to fig­ures and di­a­grams in both eps and pdf for­mats. It is suit­able both for la­belling a new di­a­gram and for re­la­belling an ex­ist­ing di­a­gram. 

\usepackage{graphicx} %The pack­age builds upon the graph­ics pack­age, pro­vid­ing a key-value in­ter­face for op­tional ar­gu­ments to the \in­clude­graph­ics com­mand. This in­ter­face pro­vides fa­cil­i­ties that go far be­yond what the graph­ics pack­age of­fers on its own. 

\usepackage[all]{xy} %A pack­age for type­set­ting a va­ri­ety of graphs and di­a­grams with TeX

\usepackage{enumitem} %This pack­age pro­vides user con­trol over the lay­out of the three ba­sic list en­vi­ron­ments: enu­mer­ate, item­ize and de­scrip­tion. It su­per­sedes both enu­mer­ate and md­wlist (pro­vid­ing well-struc­tured re­place­ments for all their fun­tion­al­ity), and in ad­di­tion pro­vides func­tions to com­pute the lay­out of la­bels, and to ‘clone’ the stan­dard en­vi­ron­ments, to cre­ate new en­vi­ron­ments with coun­ters of their own. 
\usepackage[T5,T1]{fontenc} %The pack­age al­lows the user to se­lect font en­cod­ings, and for each en­cod­ing pro­vides an in­ter­face to ‘font-en­cod­ing-spe­cific’ com­mands for each font. Its most pow­er­ful ef­fect is to en­able hy­phen­ation to op­er­ate on texts con­tain­ing any char­ac­ter in the font. 

\usepackage[utf8]{inputenc} %The pack­age trans­lates var­i­ous stan­dard and other in­put en­cod­ings into a ‘LaTeX in­ter­nal lan­guage’. The in­ter­nal lan­guage is ex­pressed en­tirely in TeX's base en­cod­ing (stan­dard ASCII print­able char­ac­ters, car­riage con­trol to­kens and TeX con­trol se­quences, the lat­ter mostly de­fined by LaTeX). 

\usepackage{lmodern} % The Latin Modern fam­ily of fonts con­sists of 72 text fonts and 20 math­e­mat­ics fonts, and is based on the Com­puter Modern fonts re­leased into pub­lic do­main by AMS (copy­right © 1997 AMS). The lm font set con­tains a lot of ad­di­tional char­ac­ters, mainly ac­cented ones, but not ex­clu­sively. There is one set of fonts, avail­able both in Adobe Type 1 for­mat (*.pfb) and in OpenType for­mat (*.otf). There are five sets of TeX Font Met­ric files, cor­re­spond­ing to: Cork en­cod­ing (cork-*.tfm); QX en­cod­ing (qx-*.tfm); TeX’n’ANSI aka LY1 en­cod­ing (tex­nansi-*.tfm); T5 (Viet­namese) en­cod­ing (t5-*.tfm); and Text Com­pan­ion for EC fonts aka TS1 (ts1-*.tfm). 

\usepackage[vietnam,english]{babel} %The pack­age man­ages cul­tur­ally-de­ter­mined ty­po­graph­i­cal (and other) rules, and hy­phen­ation pat­terns for a wide range of lan­guages. A doc­u­ment may se­lect a sin­gle lan­guage to be sup­ported, or it may se­lect sev­eral, in which case the doc­u­ment may switch from one lan­guage to an­other in a va­ri­ety of ways. 

\DeclareTextSymbolDefault{\OHORN}{T5}
\DeclareTextSymbolDefault{\UHORN}{T5}
\DeclareTextSymbolDefault{\ohorn}{T5}
\DeclareTextSymbolDefault{\uhorn}{T5}
\DeclareTextSymbolDefault{\ecircumflex}{T5}
\usepackage{amsfonts} %An ex­tended set of fonts for use in math­e­mat­ics, in­clud­ing: ex­tra math­e­mat­i­cal sym­bols; black­board bold let­ters (up­per­case only); frak­tur let­ters; sub­script sizes of bold math italic and bold Greek let­ters; sub­script sizes of large sym­bols such as sum and prod­uct; added sizes of the Com­puter Modern small caps font; cyril­lic fonts (from the Univer­sity of Wash­ing­ton); Euler math­e­mat­i­cal fonts. All fonts are pro­vided as Adobe Type 1 files, and all ex­cept the Euler fonts are pro­vided as METAFONT source. 

\usepackage{bm} %The bm pack­age de­fines a com­mand \bm which makes its ar­gu­ment bold. The ar­gu­ment may be any maths ob­ject from a sin­gle sym­bol to an ex­pres­sion. This is closely re­lated to the spec­i­fi­ca­tion of the \boldsym­bol com­mand in AMS-LaTeX, but \bm is rather more care­ful in the way it does things. 

\usepackage{mathrsfs} %Sup­port use of the Raph Smith’s For­mal Script font in math­e­mat­ics. Pro­vides a \math­scr com­mand, rather than over­writ­ing the stan­dard \math­cal com­mand, as in calrsfs.

\usepackage[mathcal]{euscript}

\usepackage{multirow} % The pack­age has a lot of flex­i­bil­ity, in­clud­ing an op­tion for spec­i­fy­ing an en­try at the “nat­u­ral” width of its text. 

\usepackage{array} %  An ex­tended im­ple­men­ta­tion of the ar­ray and tab­u­lar en­vi­ron­ments which ex­tends the op­tions for col­umn for­mats, and pro­vides "pro­grammable" for­mat spec­i­fi­ca­tions.  This pack­age is part of the la­tex-tools bun­dle in the LaTeX re­quired dis­tri­bu­tion. 

\newcolumntype{L}[1]{>{\raggedright\let\newline\\\arraybackslash\hspace{0pt}}m{#1}} %left-align
\newcolumntype{C}[1]{>{\centering\let\newline\\\arraybackslash\hspace{0pt}}m{#1}} %justify-align
\newcolumntype{R}[1]{>{\raggedleft\let\newline\\\arraybackslash\hspace{0pt}}m{#1}} %right-align
\newcolumntype{M}[1]{*1{>{\centering\arraybackslash}m{#1}}} %vertical-align, need to add @{}m{0pt}@{} at the end!
\usepackage{fancybox} % Pro­vides vari­ants of \fbox: \shad­ow­box, \dou­ble­box, \oval­box, \Oval­box, with help­ful tools for us­ing box macros and flex­i­ble ver­ba­tim macros. You can box math­e­mat­ics, floats, cen­ter, flush­left, and flushright, lists, and pages. 

\usepackage{wrapfig} % Al­lows fig­ures or ta­bles to have text wrapped around them. Does not work in com­bi­na­tion with list en­vi­ron­ments, but can be used in a par­box or mini­page, and in twocol­umn for­mat. Sup­ports the float pack­age. 

\usepackage{pgf,tikz} %  PGF is a macro pack­age for cre­at­ing graph­ics. It is plat­form- and for­mat-in­de­pen­dent and works to­gether with the most im­por­tant TeX back­end drivers, in­clud­ing pdfTeX and dvips. It comes with a user-friendly syn­tax layer called TikZ.  Its us­age is sim­i­lar to pstricks and the stan­dard pic­ture en­vi­ron­ment. PGF works with plain (pdf-)TeX, (pdf-)LaTeX, and ConTeXt. Un­like pstricks, it can pro­duce ei­ther PostScript or PDF out­put. 

\usetikzlibrary{arrows} % Arrow tip library

\usetikzlibrary{automata} % Automata Drawing Library, accessed by \usetikzlibrary{automata}, and is used for drawing "finite state automata and Turing Machines". To draw these graphs, each node, its name and relative position is defined, as well as the types of path between each.

\usetikzlibrary{backgrounds} % Background Library, accessed by \usetikzlibrary{backgrounds}, and "defines background for pictures". To use this in a Tikzpicture, an option is passed, e.g. \begin{tikzpicture}[show background rectangle], with a background rectangle style defined before the picture. (e.g. \tikzstyle{background rectangle}=[define background rectangle here]

\usetikzlibrary{calc} % Calc library, accessed throgh \usetikzlibrary{calc} to make complex coordinate caculations.

\usetikzlibrary{calendar} % Calendar Library, accessed through \usetikzlibrary{calendar}. This library is used to display calendars (I guess it's a Ronseal thing). You define a calendar as \calendar[display options and date options](Name (optional)).

\usetikzlibrary{er} % Entity Relationship Diagram Library, accessed by \usetikzlibrary{er}, as in the automata drawing library, each node is defined, as is each edge between each node, as well as any attributes. As a note of warning, underlining should be used for attributes, but this is not used as it is both ugly and difficult to implement. Italics are used instead.

\usetikzlibrary{intersections} % Intersection library, accessed through \usetikzlibrary{intersections}, to calculate intersections of paths. 

\usetikzlibrary{mindmap} % Mind map library, accessed through \usetikzlibrary{mindmap}. 

\usetikzlibrary{matrix} % Matrix Library, accessed through \usetikzlibrary{matrix}. Matrices are defined in the same way as in maths mode, however, each item in the matrix as assigned a value as a node, starting from 1. Each node can then be identified and manipulated. Delimiters can also be selected in the matrix options and can be "any delimiter that is acceptable to TeX’s \left command".

\usetikzlibrary{folding} % Paper folding library \usetikzlibrary{folding}.

\usetikzlibrary{patterns} % Pattern Library \usetikzlibrary{patterns}. This package "defines patterns for filling areas". In the documentation, each pattern is named and an example given.

\usetikzlibrary{plothandlers} % Plot handler Library, accessed through \usetikzlibrary{plothandlers}. TikZ loads this library automatically. Each point is defined (as a node) for the plot and the each point has a curve placed 

\usetikzlibrary{plotmarks} % Plot Mark Library, accessed through \usetikzlibrary{plotmarks} is used to define additional styles for plots as used above. Each point is defined as \pgfuseplotmark{Plot description}.

\usetikzlibrary{shapes} % Shape library, used to define shapes other than rectangle, circle and co-ordinate. Accessed through either \usetikzlibrary{shapes} or \usetikzlibrary{shapes.shape type}. The following additional types are available: geometric shapes, either named shapes (star, diamond, etc.) or polygons of specified side numbers; symbol shapes, e.g. "forbidden sign" as used in No Smoking signs; "multipart" shapes, with "multiple (text) parts"; and finally, "misc" shapes which "do not fit in the previous categories", such as strikethrough crosses.

\usetikzlibrary{topaths} % To Path library, accessed through \usetikzlibrary{topaths}. This library is used to define paths between two points, and is loaded automatically. Additionally, it can take the form of curved lines between two shapes or as a loop back to a node.

\usetikzlibrary{trees} % Tree library, accessed through \usetikzlibrary{trees}. Each point on the tree is defined as a node, with children, and each child can have its own children. The tree's direction can also be specified, as well as the angle at which children emerge, however, when left to its own devices, the results are acceptable.

\usetikzlibrary{decorations.pathreplacing,decorations.pathmorphing}

\usetikzlibrary{positioning,spy,shapes.geometric,fadings}
\usepackage{xcolor} % The pack­age starts from the ba­sic fa­cil­i­ties of the color pack­age, and pro­vides easy driver-in­de­pen­dent ac­cess to sev­eral kinds of color tints, shades, tones, and mixes of ar­bi­trary col­ors. It al­lows a user to se­lect a doc­u­ment-wide tar­get color model and of­fers com­plete tools for con­ver­sion be­tween eight color mod­els. Ad­di­tion­ally, there is a com­mand for al­ter­nat­ing row col­ors plus re­peated non-aligned ma­te­rial (like hor­i­zon­tal lines) in ta­bles. Colors can be mixed like \color{red!30!green!40!blue}. 

\usepackage{colortbl} % The pack­age al­lows rows and columns to be coloured, and even in­di­vid­ual cells.

\definecolor{Gray}{gray}{0.85}
\definecolor{LightCyan}{rgb}{0.88,1,1}
\definecolor{zzttqq}{rgb}{0.6,0.2,0.0}
\definecolor{qqqqff}{rgb}{0.0,0.0,1.0}
\definecolor{cqcqcq}{rgb}{0.7529411764705882,0.7529411764705882,0.7529411764705882} 
\makeatletter
\def\thm@space@setup{\thm@preskip=0pt
\thm@postskip=0pt}
\makeatother
\newtheoremstyle{newstyle}      
{} %Aboveskip 
{} %Below skip
{\itshape} %Body font e.g.\mdseries,\bfseries,\scshape,\itshape
{} %Indent
{\bfseries} %Head font e.g.\bfseries,\scshape,\itshape
{.} %Punctuation afer theorem header
{ } %Space after theorem header
{} %Heading

\newtheorem{thm}{Theorem}[section] % Standard theorem environment
\newtheorem{lem}[thm]{Lemma} % Lemma environment with numbering
\newtheorem{prop}[thm]{Proposition}
\newtheorem{cor}[thm]{Corollary} 

\theoremstyle{definition}
\newtheorem{defi}[thm]{Definition}
\newtheorem{depo}[thm]{Definition-Proposition}

\newtheorem{exam}[thm]{Example}

\newtheorem{quest}{Question}
\newtheorem{rem}[thm]{Remark}

\newtheorem*{question*}{Question}
\newtheorem*{theorem*}{Theorem}
\newtheorem*{property*}{Property}
\newtheorem*{conj*}{Conjecture}

\newtheorem*{pseudo*}{Pseudo-hyperresolution}

\newcounter{nmdthmcnt}
\newenvironment{namedthm}[2][]{\addtocounter{nmdthmcnt}{1}\theoremstyle{theorem}
	\newtheorem*{nmdthm\roman{nmdthmcnt}}{Theorem #2}%
	\begin{nmdthm\roman{nmdthmcnt}}[#1]}{\end{nmdthm\roman{nmdthmcnt}}}

\newenvironment{namedpro}[2][]{\addtocounter{nmdthmcnt}{1}\theoremstyle{theorem}
	\newtheorem*{nmdthm\roman{nmdthmcnt}}{Proposition #2}%
	\begin{nmdthm\roman{nmdthmcnt}}[#1]}{\end{nmdthm\roman{nmdthmcnt}}}

\def\ntc{\textbf{Nguy\~{\^{e}}n} Th\'{\^{e}} C\uhorn\`\ohorn ng}

\def\phr{pseudo-hyperresolution }

\def\Phr{\uppercase{p}seudo-hyperresolution }

\def\cps{infinite complex projective space }

\def\eapg{elementary abelian $p-$group}

\def\ea2g{elementary abelian $2-$group}

\def\minj{the minimal injective resolution }

\def\a2{\mathcal A_2} %mod 2 Steenrod algebra
\def\u{\mathcal U} % the category of unstable modules

\def\fk{\mathscr F_{\Bbbk}}
\def\f2{\mathbb F_2}
\def\fp{\mathbb F_p}
\def\v{\mathcal V}

\def\h{\textup{H}}
\def\n{\mathbb N}

\def\z{\mathbb Z}

\def\pc{\mathbb CP}

\def\p{\mathcal P}

\def\restrict#1{\raise-.5ex\hbox{\ensuremath|}_{#1}}

\newcommand{\gln}[2]{\mathrm{GL}_{#1}\left(#2\right)}
\newcommand{\Gln}[2]{\mathrm{GL}_{#1,#2}}
\newcommand{\pb}[2]{\mathscr{G}^{#1}\left(#2\right)}
\newcommand{\tphi}[1]{\tilde{\Phi}\left(#1\right)}
\newcommand{\moc}[1]{\left(#1\right)}
\newcommand{\mocv}[1]{\left[#1\right]}
\newcommand{\mocn}[1]{\left\{#1\right\}}

\newcommand{\ali}[1]{\begin{align*}#1\end{align*}}
\newcommand{\alilabel}[2]{\begin{align}\label{#1}#2\end{align}}
\newcommand{\exg}[3]{\mathrm{Ext}_{\u}^{#1}\left(\Sigma^{#2}\f2,\Sigma^{#3}\f2\right)}
\newcommand{\psr}[2]{\mathscr{P}\left(#1,#2\right)}
\newcommand{\hml}[2]{\mathrm{HML}^{*}\left(#1,#2\right)}
\newcommand{\Hml}[3]{\mathrm{HML}^{#1}\left(#2,#3\right)}

\newcommand{\cke}[1]{\mathrm{Coker}\left(#1\right)}

\newcommand{\ext}[4]{\mathrm{Ext}_{#1}^{#2}\left(#3,#4\right)}

\newcommand{\ho}[3]{\mathrm{Hom}_{#1}\left(#2,#3\right)}
\newcommand{\map}[3]{\mathrm{Map}_{#1}\left(#2,#3\right)}

\newcommand{\im}[1]{\mathrm{Im}\left(#1\right)}
\newcommand{\ke}[1]{\mathrm{Ker}\left(#1\right)}

\newcommand{\dsum}[3]{\underset{#1}{\overset{#2}{\bigoplus}}#3}

\newcommand{\aug}[2]{\mathrm{Aug}_{#1}\left(#2\right)}

\newcommand{\morp}[2]{#1\left(#2\right)}

\newcommand{\boundellipse}[3]% center, xdim, ydim
{(#1) ellipse (#2 and #3)
}

\newcommand{\suite}[2]{#1_{1},#1_{2},\ldots,#1_{#2}}

\newcommand\norm[1]{\left\lVert#1\right\rVert}
\numberwithin{equation}{section} % number equations by sections
\title{The Pseudo-hyperresolution and Applications}
\author{\ntc}

\title{The Pseudo-hyperresolution and Applications}
\author{\ntc}
\begin{document}
\maketitle
\begin{abstract}
Homological algebra techniques can be found in almost all modern areas of mathematics. Many interesting problems in mathematics can be formulated, computed, or can find their equivalence in terms of Ext-groups. For instance, important (co)homology theories, such as the Mac Lane cohomology for rings or the Hochschild and cyclic homology of commutative algebras can be defined as Ext-groups in suitable functor categories; homotopical invariants can also gain information from homological data with the help of the unstable Adams spectral sequence, whose input takes the form of Ext-groups in the category of unstable modules over the Steenrod algebra. Therefore, the constructions of explicit injective (projective) resolutions in an abelian category is of great importance. In this article, we introduce a new method, called  \textit{Pseudo-hyperresolution}, to study such constructions. This method originates in the category of unstable modules, and aims at building explicit resolutions for the reduced singular cohomology of spheres. In particular, for all integers $n\geq 0$, we can describe a large range of the minimal injective resolution of the sphere $S^{n}$ based on the Bockstein operation of the Steenrod algebra. Moreover, many classical constructions in algebraic topology, such as the algebraic EHP sequence or the Lambda algebra can be recovered using the \Phr method. A particular connection between spheres and the \cps is also established.  Despite its origin, Pseudo-hyperresolution generalizes to all abelian categories. In particular, many explicit resolutions of classical strict polynomial functors can be reunified in view of Pseudo-hyperresolution. As a consequence, we recover the global dimension of the category of homogeneous strict polynomial functors of finite degree as well as the Mac Lane cohomology of finite fields.
\end{abstract}

\section{Introduction}
Singular cohomology plays an important role in algebraic topology and a great deal of effort was put into developing a suitable framework to study this cohomology theory. In the process, Steenrod \cite{Ste62} defined the notion of a \textit{stable cohomology operation}, which led to the construction of the algebra that bears his name. The singular cohomology of a topological space admits a natural action of the Steenrod algebra satisfying an extra axiom that Steenrod called the \textit{instability} condition. He then introduced the notion of an \textit{unstable module} that captures this particular property. It was in the language of unstable modules that many important results in homotopy theory found their solutions, \textit{e.g.}, the Segal conjecture \cite{Car83}, the Sullivan conjecture \cite{Mil84}, the Serre conjecture \cite{LS86}, the Kuhn realization conjecture \cite{Kuh95}, \textit{etc}. As their name suggests, unstable modules are an important tool in studying unstable homotopy theory, and in particular, to study the homotopy groups of topological spaces. Among all the interesting spaces in topology, the
spheres are the most fundamental and important. From the categorical point of view, we should not just focus on the
objects themselves, but the maps between them as well. To this end, in algebraic topology, we seek to classify the set
of continuous maps between spheres up to continuous deformation. The importance of this classification is not just
aesthetical, but also from the fact that it is connected to many other areas in mathematics, such as geometric topology,
algebra, and algebraic geometry. For example, the groups of differential structures on spheres is determined by the stable homotopy groups of spheres (see, \textit{e.g.}, \cite{KM63}). Another example is the theory of topological modular forms (see, \textit{e.g.}, \cite{Hop02}), which relates certain parts of the stable homotopy groups of spheres to the moduli stack of elliptic curves. However, the homotopy groups of spheres are a hugely intractable object even though their cohomology groups are elementary. One of the most powerful tool in studying the homotopy groups of spheres is the unstable Adams spectral sequence - as it was introduced by Massey and Peterson in \cite{MP67}, generalized by Bousfield and Curtis in \cite{BC70}, and generalized further by Bousfield and Kan in \cite{BK72} - which
passes from homological information to homotopical information. In fact, it reduces the computations of the homotopy
groups of spheres to that of certain Ext-groups in the category $\u$ of unstable modules:
\ali{E_{2}^{s,t}\moc{S^{n}}=\ext{\u}{s}{\morp{\tilde{\mathrm{H}}^{*}}{S^{n};\z/p}}{\morp{\tilde{\mathrm{H}}^{*}}{S^{t};\z/p}}\Longrightarrow \morp{\pi_{t-s}}{S^{n}}^{\wedge}_{p}.}
Therefore, the understanding of the minimal injective resolution of $\morp{\tilde{\mathrm{H}}^{*}}{S^{t};\z/p}$ is of vital importance, and is one of the primary goals of this article. Our approach to this problem relies on the simple structure of $\morp{\tilde{\mathrm{H}}^{*}}{S^{t};\z/p}$: it is an $\mathbb{N}-$graded $\fp-$vector space, concentrated in degree $t$ and isomorphic to $\fp$ in that degree, and therefore is denoted by $\Sigma^{t}\fp$. Moreover, the tensor product with $\Sigma^{t}\fp$ defines an exact functor $\Sigma^{t}:\u\to\u$, called the $t-$th suspension functor (we simply write $\Sigma$ when $t=1$ and call it the suspension functor). Now, it is well-known that $\Sigma^{t}\fp$ is an injective unstable module for $t=0$ and $t=1$. This suggests that the construction of \minj of $\Sigma^{t}\fp$ should be carried out by induction on $t$ and this is what we are going to do. In fact, as the suspension functor $\Sigma$ is exact and $\Sigma\moc{\Sigma^{t}}=\Sigma^{t+1}$, then by applying this functor to an injective resolution of $\Sigma^{t}\fp$ we obtain an acyclic cochain complex that admits $\Sigma^{t+1}\fp$ as its only nontrivial cohomology. Even though this resulting complex is no longer an injective resolution, each of its terms is the suspension of an injective unstable module. On the other hand, such an unstable module has a simple injective resolution induced by the \textit{Mahowald exact sequences}. It turns out that these resolutions suffice to construct an explicit injective resolution for $\Sigma^{t+1}\fp$ with the help of the \textit{\Phr}method, which is the solution to the following general question.
\begin{quest}\label{generalqu}
Let $\mathcal{C}$ be an abelian category with enough injective objects and let $M\in\mathcal{C}$ be the only nontrivial cohomology of a certain cochain complex $\moc{A^{k},\partial^{k}}_{k\geq 0}$. Suppose that $A^{k}$ admits an explicit injective resolution in $\mathcal{C}$ for all $k\geq 0$. Then, is it possible to construct an explicit injective resolution for $M$, based on these given resolutions of all $A^{k}$? 
\end{quest}
The \phr method, as it will be introduced in Section \ref{Phr}, plays an essential role in the present paper. So let recall how it is formulated. An easy way to describe the \phr method is to use the Poincar\'{e} power series of (co)chain complexes. That is, given a (co)chain complex $\moc{A^{k},\partial^{k}}_{k\geq 0}$, we define its power series $\morp{\mathrm{P}}{A^{k},\partial^{k}}$ as
\ali{\morp{\mathrm{P}}{A^{k},\partial^{k}}:=\sum_{k=0}^{\infty}t^{k}\mocv{A^{k}},}
where $\mocv{A^{k}}$ denotes the corresponding class of $A^{k}$ in the Grothendieck group of $\mathcal{C}$. 
The answer to Question \ref{generalqu} now goes as follows. 
\begin{pseudo*}
Denote by $\moc{I^{k,\bullet},\partial^{k,\bullet}}$ the given injective resolution of $A^{k}$ for all integers $k\geq 0$. Then, there exists an injective resolution $\moc{J^{\bullet},\delta^{\bullet}}$ of $M$ such that
\ali{\morp{\mathrm{P}}{J^{\bullet},\delta^{\bullet}}=\sum_{k=0}^{\infty}t^{k}\morp{\mathrm{P}}{I^{k,\bullet},\partial^{k,\bullet}}.}
\end{pseudo*}
Note that, when the length of the (co)chain complex in play is finite, the resulting pseudo-hyperresolution is well-known to algebraists: this resolution is identical to the one obtained using the successive cone technique. However, for infinite-length complexes, algebraists have to use the mapping telescope technique, which yields more complicated resolutions than those induced by the pseudo-hyperresolution method.

Before explaining why the \phr method is useful in studying unstable modules, let recall some basic facts about these objects. Note that, most of the methods we use throughout this article works for all prime characteristics. Some particular methods are only treated in characteristic $2$, but there is no difficulty extending to all prime characteristics. However, special care is required in dealing with odd prime characteristics. Therefore, they will be studied separately in a subsequent article. Moreover, when it comes to concrete examples and recalls, we always restrict the attention to characteristic $2$ for ease of exposition. 

The mod $2$ Steenrod algebra is the quotient of the free associative unital graded $\f2-$algebra generated by the symbols $Sq^{k}$ of degree $k\geq 0$, subject to the Adem relations. An unstable module $M$ is an $\mathbb{N}-$graded module over the Steenrod algebra such that for all elements $x\in M$ of degree $n$, then $Sq^{k}x=0$ for all integers $k>n$. This particular property of the Steenrod action is verified by the singular cohomology of a space but not that of a spectrum, whence the appellation \textit{unstable}. There are two typical families of injective unstable modules. The first one consists of the singular cohomology of \ea2gs and the second one consists of the injective envelope $J(n)$ (also known as the $n-$th Brown-Gitler module) of the reduced singular cohomology of the sphere $S^{n}$ for all integers $n\geq 0$. The discovery of the first family is not trivial. It is the essential key to the Segal and Sullivan conjectures, two
of the great highlights in algebraic topology in the late 1980’s and 1990’s, which led to the establishment of a whole new
area within the subject. In fact, Carlsson \cite{Car83} observed the injectivity of the cohomology $\h^{*}\moc{B\z/p;\fp}$ of the classifying space of $\z/p$ for $p=2$. Miller \cite{Mil84} then extended this result to all prime $p$. In \cite{LZ86}, Lannes and Zarati made a spectacular contribution to the theory of unstable modules by showing that the tensor product  $\h^{*}\moc{BV;\fp}\otimes J(n)$ remains injective in $\u$ for all \eapg{s} $V$ and all integers $n\geq 0$. Finally, the full characterization of injective unstable modules was achieved in \cite{LS89}, which amounts to saying that all injective unstable module is isomorphic to a direct sum of tensor products of the form  $\h^{*}\moc{BV;\fp}\otimes J(n)$. Now, let get back to the second family of injective unstable modules. The Brown-Gitler module $J(n)$ is the representing object of the functor $M\mapsto \ho{\f2}{M^{n}}{\f2}$ from $\u$ to the category of $\f2-$vector spaces. Hence, they form a system of injective cogenerators for the category $\u$. Moreover, Brown-Gitler modules are connected by the short Mahowald exact sequences, which amounts to saying that the suspension $\Sigma J(2k)$ is again an injective unstable module as it is isomorphic to $J(2k+1)$, whereas $\Sigma J(2k+1)$ is of injective dimension $1$ and its minimal injective resolution is given by $J(2k+2)\twoheadrightarrow J(k+1)$, which is induced by the Steenrod operation $Sq^{k+1}$. It follows from the characterization of injective unstable modules that the suspension of an injective object of $\u$ is of injective dimension at most $1$ and admits an explicit injective resolution induced by the Mahowald short exact sequences. Therefore, in view of pseudo-hyperresolution, if an unstable module admits an explicit injective resolution, then so is its suspension. For instance, the reduced singular cohomology of spheres and that of the \cps are such objects. Let discuss about the former first. 

As $\morp{\tilde{\h}^{*}}{S^{1};\f2}\cong J(1)$, and $\Sigma\morp{\tilde{\h}^{*}}{S^{n};\f2}\cong \morp{\tilde{\h}^{*}}{S^{n+1};\f2}\cong \Sigma^{n+1}\f2$ for all integers $n\geq 1$, then the \phr method allows to construct for each integer $n\geq 1$ an explicit injective resolution of $\Sigma^{n}\f2$. In general, the \phr method fails to provide information about differentials. However, in this particular case, we have enough control on the resulting resolution of $\Sigma^{n}\f2$ allowing to obtain the minimal one in a large range. This is carried out with the help of the Bockstein operation of the Steenrod algebra. In fact, recall that a morphism between two Brown-Gitler modules $J(n)$ and $J(m)$ is determined by a Steenrod operation $\theta$ of degree $n-m$, and we denote this morphism by $\bullet \theta$. As the square power of $Sq^{1}$ is trivial, then the sequence $\moc{J(n),\bullet Sq^{1}}_{n\geq 1}$ is a complex, which is very close to being exact.
\begin{namedpro}{\ref{bsss}}
	Let $n\geq 1$ be an integer. Then, the  sequence
	\begin{equation*}
	J(n+2)\oplus \mathscr{J}\left( n+2\right)\xrightarrow{\left(\begin{smallmatrix}
		\bullet Sq^{1}&0\\\bullet Sq^{\frac{n+2}{2}}&0
		\end{smallmatrix}\right)}
	J(n+1)\oplus \mathscr{J}\left( n+1\right)\xrightarrow{\left(\begin{smallmatrix}
		\bullet Sq^{1}&0\\ \bullet Sq^{\frac{n+1}{2}}&0
		\end{smallmatrix}\right)}
	J(n)\oplus \mathscr{J}\left( n\right)
	\end{equation*}
	is exact, where 
	$$\mathscr{J}\left( n\right)=\left\{\begin{array}{cl}
	J(2k)&\text{ if }n+1=4k,\\
	0 &\text{ otherwise.}
	\end{array}
	\right.$$
\end{namedpro} 
It turns out that this particular long exact sequence agrees with the minimal injective resolution of $\Sigma^{t}\f2$ in a large area:
     	 \begin{namedthm}{\ref{ressphere}}
     	 	Let $n,s,t\geq 0$ be three  integers such that $s> [n/2]$, where $[-]$ denotes the integral part of a number. Then, we have
     	 	\begin{equation*}
     	 	\ext{\u}{s}{\Sigma^{n}\f2}{\Sigma^{t}\f2}\cong\left\{\begin{array}{ll}
     	 	\f2&\text{ if }n=t-s,\\
     	 	\f2&\text{ if }t-s-1\equiv 0 (4)\text{ and }t-s-1=2n,\\
     	 	0&\text{ otherwise.}
     	 	\end{array}\right.
     	 	\end{equation*}
     	 	In other words, in this case, the $s-$th term of the minimal injective resolution of $\Sigma^{t}\f2$ is isomorphic to 
     	 	\ali{\morp{J}{t-s}\oplus \morp{\mathscr{J}}{t-s}.}
     	 \end{namedthm}
An alternative way to compute $\ext{\u}{s}{\Sigma^{n}\f2}{\Sigma^{t}\f2}$ relies on the fact that $\ho{\u}{\Sigma^{n}\f2}{\morp{J}{k}}$ is isomorphic to $\f2$ if $k=n$ and is trivial otherwise. Therefore, $\ext{\u}{s}{\Sigma^{n}\f2}{\Sigma^{t}\f2}$ is isomorphic to $\f2^{\oplus d}$, where $d$ is the number of copies of direct summands of the form $J(n)$ in the $s-$th term of the minimal injective resolution of $\Sigma^{t}\f2$. Because the \phr method passes information from the minimal injective resolution of $\Sigma^{t}\f2$ to that of the minimal injective resolution of $\Sigma^{t+1}\f2$, then by using Mahowald short exact sequences - which is the algebraic analogue of the James fibrations - we recover the algebraic EHP sequence.
\begin{namedthm}{\ref{EHP}}
	There exists a long exact sequence
	\begin{equation*}
	E_{2}^{s-2,t}(S^{2n+1})\xrightarrow{P}E_{2}^{s,t}(S^{n})\xrightarrow{E}E_{2}^{s,t+1}(S^{n+1})\xrightarrow{H}E_{2}^{s-1,t}(S^{2n+1})
	\end{equation*}
    for all integers $n\geq 0$ and $s\geq 2$, where $E_{2}^{s,t}(S^{n})$ stands for $\ext{\u}{s}{\Sigma^{n}\f2}{\Sigma^{t}\f2}$. 
\end{namedthm}

In fact, the algebraic EHP sequence exists in a much more general form. Recall that the suspension functor $\Sigma$ admits a left adjoint, denoted by $\Omega$. Having observed that the left-derived functors of $\Omega$, denoted by $\Omega_{s}$ for all integers $s\geq 0$, are trivial in homological degrees greater than $1$, Bousfield showed that there existed an exact sequence $$\ext{\u}{s-2}{\Omega_{1}M}{N}\xrightarrow{}\ext{\u}{s}{\Omega M}{N}\xrightarrow{}\ext{\u}{s}{M}{\Sigma^{}N}\xrightarrow{}\ext{\u}{s-1}{\Omega_{1}M}{N}$$
for all integers $s\geq 2$ and all unstable modules $M,N$. His proof might be well-known to experts, but is not available in the literature. (As pointed out to the author by Hans-Werner Henn, only a compact explanation can be found in \cite{Mil84}.) Following Singer's instruction \cite{Wil16}, we also recover Bousfield's proof. Thanks to the generous permission of Bousfield, we give a detailed account of his approach in Section \ref{BS}. 

Another application of the \phr method that we discuss in this article is the construction of the Lambda algebra. In stead of studying the minimal injective resolution, we construct for each integer $t\geq 0$ a particular injective resolution of $\Sigma^{t}\f2$. These resolutions fit together into a direct systems of which the limit is endowed with the structure of a bigraded differential algebra, which is isomorphic to the Lambda algebra. This construction can then be considered as a \textit{special injective resolution of} $\Sigma^{\infty}\f2$. 
\begin{namedpro}{\ref{lambdaalgebra}}
There exists a bigraded differential algebra $\moc{\Lambda,d}$ such that 
\ali{\h^{r,s}\moc{\Lambda,d}\cong \varinjlim\limits_{n}\ext{\u}{r}{\Sigma^{n}\f2}{\Sigma^{n+s}\f2}.}
\end{namedpro}

We now turn to the case of the reduced singular cohomology of the infinite complex projective space. Recall that it is isomorphic to the augmentation ideal of the polynomial algebra on one variable of degree $2$. As the cohomology $\h^{*}\moc{B\z/2;\f2}$ of the classifying space of $\z/2$ is isomorphic to the polynomial algebra on one variable of degree $1$, we obtain the following short exact sequence.
\ali{0\to \Sigma^{t}\tilde{\h}^{*}\moc{\pc^{\infty};\f2}\to \Sigma^{t}\tilde{\h}\moc{B\z/2;\f2}\to \Sigma^{t+1} \h^{*}\moc{\pc^{\infty};\f2}\to 0.} 
As a result, $\Sigma^{t}\tilde{\h}^{*}\moc{\pc^{\infty};\f2}$ is the only nontrivial cohomology of the acyclic cochain complex
\begin{align}\label{compsph1}
\Sigma^{t}\bar{\h}\to \Sigma^{t+1} \h\to\cdots\to\Sigma^{t+k}\h\to\Sigma^{t+k+1}\h\to\cdots,
\end{align}
where $\bar{\h}$ stands for $\tilde{\h}\moc{B\z/2;\f2}$ and $\h$ stands for $\h^{*}\moc{B\z/2;\f2}$. As a consequence, we obtain the spectral sequence
\begin{align}\label{sscomp}
E_{2}^{q,k}:=\left\{ \begin{array}{cl}
\ext{\u}{q}{\Sigma^{n}\f2}{\Sigma^{t+k}\h}&\textup{ if }k\geq 1,\\
\ext{\u}{q}{\Sigma^{n}\f2}{\Sigma^{t}\bar{\h}}&\textup{ if }k=0.
\end{array} \right.\Longrightarrow \ext{\u}{q+k}{\Sigma^{n}\f2}{\Sigma^{t}\tilde{\h}^{*}\moc{\pc^{\infty};\f2}}.
\end{align}
Moreover, because of the characterization of injective unstable modules (see, \textit{e.g.}, \cite{LS89}), then the tensor product of an explicit injective resolution of $\Sigma^{t+k}\f2$ with $\h$ (or $\bar{\h}$) yields an explicit injective resolution of $\Sigma^{t+k}\h$ (or $\Sigma^{t+k}\bar{\h}$). Therefore, in view of \phr method, $\Sigma^{t}\tilde{\h}^{*}\moc{\pc^{\infty};\f2}$ admits an explicit injective resolution. By inspecting this resolution, we show that the spectral sequence \eqref{sscomp} collapses at the $E_{2}-$term, giving rise to the isomorphism
\begin{align}\label{compsph}
\ext{\u}{s}{\Sigma^{n}\f2}{\Sigma^{t}\tilde{\h}^{*}\moc{\pc^{\infty};\f2}}&\cong\ext{\u}{s}{\Sigma^{n}\f2}{\Sigma^{t}\bar{\h}}\oplus\bigoplus_{m=1}^{s} \ext{\u}{s-m}{\Sigma^{n}\f2}{\Sigma^{t+m}\h}.
\end{align}
As $\h$ is isomorphic to $\bar{\h}\oplus\f2$ and $\ho{\u}{\Sigma^{n}\f2}{\bar{H}\otimes J(k)}\cong 0$ for all integers $n,k\geq 0$, then we obtain the following $\f2-$isomorphism.
\begin{namedthm}{\ref{compsphfinal}}
For all integers $s\geq 0$, we have
\begin{align}\label{compsph2}
\ext{\u}{s}{\Sigma^{n}\f2}{\Sigma^{t}\tilde{\h}^{*}\moc{\pc^{\infty};\f2}}&\cong \bigoplus_{m=1}^{s} \ext{\u}{s-m}{\Sigma^{n}\f2}{\Sigma^{t+m}\f2}.
\end{align}
\end{namedthm}
The intrigued readers might wonder what is the interest in studying such a relation. Here is the reason. Because, on the one hand, there exists a spectral sequence
\ali{E_{2}^{s,t}:=\ext{\u}{s}{\Sigma^{n}\f2}{\Sigma^{t}\tilde{\h}^{*}\moc{\pc^{\infty};\f2}}\Longrightarrow \morp{\pi_{t-s}}{\map{*}{\pc^{\infty}}{\hat{S^{n}}}},}
where $(-)^{\wedge}$ denotes the profinite completion of a topological space, and on the other hand, the mapping space $\map{*}{\pc^{\infty}}{\hat{S^{n}}}$ is contractible (see, \textit{e.g.}, \cite{Zab87,McG96}), then the isomorphism \eqref{compsph2} is of interest as it might give new information about the groups $\ext{\u}{s}{\Sigma^{n}\f2}{\Sigma^{t}\f2}$.

Despite its origin, pseudo-hyperresolution generalizes to all abelian categories. Of particular interest is the category of strict polynomial functors, which is closely related to the category of unstable modules (see, \textit{e.g.}, \cite{Hai10,Ngu16}). The notion of a strict polynomial functor, as it is introduced in \cite{FS97}, plays a central role in the proof of Friedlander and Suslin that for a finite group scheme $G$ and a finite dimensional rational module $M$, then  $\h^{*}\left(G;\Bbbk\right)$ is an algebra of finite type and $\h^{*}\left(G;M\right)$ is a module of finite type over $\h^{*}\left(G;\Bbbk\right)$, where $\Bbbk$ is a finite field. A strict polynomial functor is a functor from the category of $\Bbbk-$vector spaces of finite dimension to the category of  $\Bbbk-$vector spaces such that for all couple of $\Bbbk-$vector spaces $(V,W)$, the structural morphism from $\ho{\Bbbk}{V}{W}$ to $\ho{\Bbbk}{F(V)}{F(W)}$ is a scheme morphism. Here, we identify  a $\Bbbk-$vector space $V$ with the affine scheme $Spec\left(S^{*}(V^{\sharp})\right)$, where $(-)^{\sharp}$ denotes the $\Bbbk-$linear dual and $S^{*}$ denotes the symmetric algebra. Let $Sch/\Bbbk$ denote the category of schemes over $\Bbbk$. Because
\ali{\ho{Sch/\Bbbk}{Spec\moc{S^{*}\moc{X^{\sharp}}}}{Spec\moc{S^{*}\moc{Y^{\sharp}}}}\cong S^{*}\moc{X^{\sharp}}\otimes Y}
for all $\Bbbk-$vector spaces $X$ and $Y$, then we define a scheme map $p:X\to Y$ to be homogeneous of degree $d$ if it belongs to $S^{d}\moc{X^{\sharp}}\otimes Y$. And then, a strict polynomial functor is defined to be homogeneous of degree $d$ if its structural morphisms are homogeneous of the same degree. Denote by $\p_{\Bbbk,d}$ the full subcategory of homogeneous strict polynomial functors of degree $d$. Then,  $\p_{\Bbbk}\cong\bigoplus_{d\geq 0}\p_{\Bbbk,d}$. In this article, we focus on the applications of the \phr method to strict polynomial functors in characteristic $2$. 

Typical examples of strict polynomial functors are given by the symmetric powers $S^{n}$, the exterior powers $\Lambda^{n}$, and the divided power $\Gamma^{n}$ for all integers $n\geq 0$. The category $\p_{\f2}$ is well-known for its computability as the injective and the projective objects are well understood. Each injective strict polynomial functor is isomorphic to a direct sum of functors of the form $S^{\lambda_{1}}\otimes S^{\lambda_{2}}\otimes\ldots\otimes S^{\lambda_{k}}$, and  each projective one is isomorphic to a direct sum of functors of the form $\Gamma^{\lambda_{1}}\otimes \Gamma^{\lambda_{2}}\otimes\ldots\otimes \Gamma^{\lambda_{k}}$, where $\moc{\suite{\lambda}{k}}\in\mathbb{N}^{\times k}$. Now, the pseudo-hyperresolution method finds its place in the category of strict polynomial functors because of the Koszul exact sequences, which connect the exterior power functors to the symmetric and the divided ones:
\ali{\xymatrixrowsep{1ex}\xymatrix{
0\ar[r]&\Lambda^{n}\ar[r]& \Lambda^{n-1}\otimes S^{1}\ar[r]&\Lambda^{2}\otimes S^{2}\ar[r]&\cdots\ar[r]& \Lambda^{1}\otimes S^{n-1}\to S^{n}\ar[r]& 0\\
0\ar[r]&\Gamma^{n}\ar[r]& \Gamma^{n-1}\otimes \Lambda^{1}\ar[r]&\Gamma^{2}\otimes \Lambda^{2}\ar[r]&\cdots\ar[r]& \Gamma^{1}\otimes \Lambda^{n-1}\to \Lambda^{n}\ar[r]& 0.
}}
As $S^{1}=\Gamma^{1}=\Lambda^{1}$, then these functors are both injective and projective. Therefore, an easy induction on $d$ using the pseudo-hyperresolution method allows to construct for each integer $d\geq 1$ an explicit (projective) injective resolution of $\Lambda^{d}$. Similarly, we can construct for each integer $d\geq 0$ an explicit injective resolution of $\Gamma^{d}$ and an explicit projective resolution of the functor $S^{d}$. These resolutions allow to show that the injective dimension of a projective strict polynomial functor of degree $d$ is bounded by $2d-2k$ and the projective dimension of an injective one is also bounded by the same number, where $k$ is the number of nontrivial $2-$adic digits of $d$. This concludes the global dimension of the category $\p_{\f2,d}$.
\begin{namedthm}{\ref{globaldimension}}
Let $d\geq 1$ be an integer and let $\sum\limits_{i=1}^{k}2^{n_{k}}$ be its $2-$adic expression. Then, the global dimension of $\p_{d}$ is $2d-2k$.
\end{namedthm}
One of the most fundamental notion in the theory of strict polynomial functors is the Frobenius twist. Let $I=\Gamma^{1}$, then we define the Frobenius twist of $I$, denoted by $I^{(1)}$, as the strict polynomial functor that associates an $\f2-$vector space $V$ to the $\f2-$vector space $V^{(1)}$, which is obtained from $V$ by base change along the Frobenius map $\f2\to\f2,x\mapsto x^{2}$. Then, we define the Frobenius twist of a functor $F$, denoted by $F^{(1)}$, as the precomposition $F\circ I^{(1)}$. The Frobenius twist is of great importance in studying strict polynomial functors for many reasons that we will recall in Section \ref{functor}. For instance, as the Frobenius twist is exact, then it induces the morphism
\alilabel{sysfrob}{\ext{\p_{\f2}}{*}{I^{(r)}}{I^{(r)}}\to\ext{\p_{\f2}}{*}{I^{(r+1)}}{I^{(r+1)}}}
for all integers $r\geq 0$. Here $I^{(r)}$ is defined recursively by $I^{(r)}=\moc{I^{(r-1)}}^{(1)}$. The colimit of  $\ext{\p_{\f2}}{*}{I^{(r)}}{I^{(r)}}$ with respect to the direct system induced by the morphisms \ref{sysfrob} is of great interest as it is isomorphic to the Mac Lane cohomology $\mathrm{HML}^{*}\moc{\f2}$ of the field $\f2$ (see, \textit{e.g.}, \cite{FS97}). It follows that the computation of $\mathrm{HML}^{*}\moc{\f2}$ can be deduced from that of $\ext{\p_{\f2}}{*}{I^{(r)}}{I^{(r)}}$ for all integers $r\geq 0$, which, in turn, can be carried out by constructing explicit injective resolutions of $I^{(r)}$. Now, recall that the Frobenius twist plays a similar role to the injective strict polynomial functors as the suspension functor does to injective unstable modules. In fact, the Frobenius twist of an injective object of $\p_{\f2}$ is no longer injective, but admits an explicit injective resolution thanks to the following exact sequence:
\ali{0\to S^{d(1)}\to S^{2d}\to S^{2d-1}\otimes S^{1}\to S^{2d-2}\otimes S^{2}\to\cdots\to S^{1}\otimes S^{2d-1}\to S^{2d}\to 0.}
Therefore, in view of \phr method, $I^{(r)}$ has an explicit injective resolution for all integers $r\geq 1$. This allows to show that 
\ali{\ext{\p_{\f2}}{k}{I^{(r)}}{I^{(r)}}\cong \left\{\begin{array}{cl}
\f2&\text{ if }2|k, \text{ and }k\leq 2^{r+1}-2,\\
0&\text{ otherwise,}
\end{array} \right.
}
which allows recovery of the Mac Lane cohomology of $\f2$:
\begin{namedthm}{\ref{mclane}}
We have
\ali{\Hml{k}{\f2}{I}\cong \left\{\begin{array}{cl}
\f2&\text{ if }2|k,\\
0&\text{ otherwise.}
\end{array} \right.
}
\end{namedthm}
The paper is organized as follows:

Section \ref{section2} recalls basic facts about the Steenrod algebra and unstable modules. We also give a brief account of the Brown-Gitler modules and show how to fit them together into the long exact sequence that we call the \textit{Bockstein sequence}.

Section \ref{graphrep} introduces the notion of a graph representation of complexes, which facilitates the presentation of our exposition. 

Section \ref{Phr} is at the heart of the present paper as it settles down the notion of a pseudo-hyperresolution. We begin this section with a simple example that motivates our study, and generalize the idea to the most general context.

Section \ref{functor} aims at providing applications of the \phr method in functor homology. It is in this section that we recover the Mac Lane cohomology of finite fields as well as the global dimension of the category of homogeneous strict polynomial functors of finite degree.

Section \ref{Lambdaalgebra} covers the construction of the Lambda algebra. It will be carried out with the help of the graph representation of complexes that we introduce in Section \ref{graphrep}.

Sections \ref{minimalres} and \ref{ressph} deal with the minimal injective resolution of the reduced singular cohomology of spheres. This will be done by studying the Bockstein sequence introduced in Section \ref{section2}.

Section \ref{ehps} provides our approach to the algebraic EHP sequence and Section \ref{BS} recalls that of Bousfield.

Section \ref{icps} studies the minimal injective resolution of the \cps and its relation to that of spheres.

\section{The Steenrod algebra and Brown-Gitler modules}\label{section2}
\noindent \textsc{\bfseries The Steenrod algebra. }The mod $2$ Steenrod algebra $\a2$ is the quotient of the free associative unital graded $\f2-$algebra generated by the symbols
$Sq^{k}$ of degree $k\geq 0$, subject to the Adem relations:
$$Sq^{i}Sq^{j}=\sum_{t=0}^{\left[\frac{i}{2}\right]}\binom{(j-t)-1}{i-2t}Sq^{i+j-t}Sq^{t}$$
for all integers $i\leq 2j$, and $Sq^{0}=1$ \cite{Ste62,Sch94}. Here, $[-]$ stands for the integral part of a number. 

An $\n-$graded  $\a2-$module $M$ is called unstable
if $Sq^{k}x=0$ for all $x\in M^{n}$ and all $k>n$. We denote by
$\u$ the category of unstable modules. Serre \cite{Ser53} introduced the notions of \textit{admissible} and \textit{excess}. A monomial 
$$Sq^{i_{1}}Sq^{i_{2}}\ldots Sq^{i_{k}}$$
is called admissible if $i_{j}\geq 2i_{j+1}$ for all $k-1\geq j\geq 1$ and $i_{k}\geq 1$; the
excess of this operation is
defined by $$e\left( Sq^{i_{1}}Sq^{i_{2}}\ldots Sq^{i_{k}}\right)=
2i_{1}-\left(\sum_{j=1}^{k}i_{j}\right).$$ The set
of admissible monomials and $Sq^{0}$ is an $\f2-$basis of $\a2$.

\noindent\textsc{\bfseries Projective unstable modules. }Recall that the functor that associates an unstable module with the $\f2-$vector space of its elements of degree $n$ is representable and we denote by $F(n)$ the representing unstable module.  Hence, $F(n)$ is freely generated by an element $\imath_{n}$ of degree $n$. Therefore, if $M$ is an unstable module, then the morphism
$$\bigoplus_{n}\bigoplus_{x\in M^{n}}F(n)\to M,\hspace{2mm}\imath_{n}\mapsto x$$
is surjective. It follows that the modules  $F(n),n\geq 0,$ form a system of projective generators of $\u$.

\noindent\textsc{\bfseries Brown-Gitler modules. }The category $\u$ also has enough injective objects. The Brown-Gitler module $J(n)$ is the unstable injective hull of the mod $2$ reduced cohomology $\Sigma^{n}\f2$ of the sphere $S^{n}$. It is a cocyclic unstable module, cogenerated by an element of degree $n$. The modules $J(n),n\geq 0,$ are
injective satisfying $\ho{\u}{M}{J(n)}\cong \left(M^{n}\right)^{\sharp}$. Here, $(-)^{\sharp}$ stands for the $\f2-$linear dual. Moreover, each $J(n)$ is $\u-$indecomposable. Each homogeneous element $x\in M^{n}$ determines a morphism $i_{x}:M\to J(n).$ Because the induced
morphism
$$i_{M}:M\to \prod_{n}\prod_{x\in M^{n}}J(n)$$
is injective, the modules $J(n),n\geq 0$, form a system of injective co-generators of $\u$.

A morphism $F(n)\to F(m)$ is determined by a Steenrod operation of degree $n-m$. Since $J(n)^{m}\cong \left(F(m)^{n}\right)^{\sharp}$, then a morphism $J(n)\to J(m)$ is
also determined by a Steenrod operation of degree $n-m$. Let
$\theta$ be a Steenrod operation of degree $n-m$, we denote by $\bullet \theta$ the morphism $J(n)\to
J(m)$ induced by $\theta$.

\noindent\textsc{\bfseries The Frobenius twist of unstable modules. }Let $|-|$ be the degree of a homogeneous element. Denote by $\mathrm{Sq}_{0}$ the operation which
associates a homogeneous element $x\in M^{n}$ of an unstable module $M$ with the element $Sq^{|x|}x$. Let $\Phi$ be the endofunctor of $\u$ which associates an
unstable module $M$ with the module $\Phi M$ concentrating in even degrees and $\left(\Phi
M\right)^{2n}=M^{n}.$ If $x\in M^{n}$, then we denote by $\Phi x$ the corresponding element in $\moc{\Phi M}^{2n}$; the action of the Steenrod algebra on $\Phi M$ is defined by $Sq^{k}\Phi x=\Phi Sq^{k/2}x$ (with the convention that $Sq^{k/2}=0$ if $k$ is not divisible by $2$). Let $\lambda_{M}$ be the morphism defined by:
\alilabel{double}{
	\lambda_{M}:\Phi M&\to M\\
	\Phi x&\mapsto \mathrm{Sq}_{0}x.\notag}
The morphism $\lambda_{M}$ is natural in $M$. There exists a right adjoint $\tilde{\Phi}$ of $\Phi$ \cite{Sch94}. Therefore,
adjoint to $\lambda_M,$ there exists a morphism $\tilde{\lambda}_M:M\to \tilde{\Phi}M,$ natural in $M$. An unstable module $M$ is reduced if $\lambda_{M}$ is injective (this is slightly different for odd prime characteristic, see, \textit{e.g.}, \cite{Sch94}). It is called nilpotent if for all $x\in M^{n},$ there exists an integer $n_{x}\geq 1$ such that $\mathrm{Sq}_{0}^{n_{x}}x=0$. By definition, there is no nontrivial morphism from a nilpotent unstable module to a reduced one.

\noindent\textsc{\bfseries The suspension functor and Mahowald's short exact sequences. }Let $\Sigma$ be the endofunctor of $\u$ defined by $(\Sigma M)^{n}=M^{n-1}$ for all unstable modules $M$ and all integers $n\geq 0$. If $x\in M^{n}$, we denote by $\Sigma x$ the corresponding element in $\Sigma M$. The action of the Steenrod algebra on $\Sigma M$ is determined by $Sq^{k}\Sigma x=\Sigma Sq^{k}x$. The suspension functor $\Sigma$ admits a right adjoint, denoted by  $\tilde{\Sigma}$. Let
$s:\Sigma\tilde{\Sigma}\to Id$ be the counit of the adjunction, and we shorten $s\moc{J(n)}$ to $s_{n}$ for all integers $n\geq 1.$
\begin{thm}[Mahowald's short exact sequences]\label{mah}
	The sequence 
	\alilabel{mahowalds}{0\to \Sigma\tilde{\Sigma}
	J(n)\xrightarrow{s_{n}}J(n)\xrightarrow{\tilde{\lambda}_{J(n)}} \tilde{\Phi}J(n)\to 0}
	is exact.
\end{thm}
Let $n\geq 0$ be an integer, we fix:
$$J\left(\frac{n}{2}\right)=\left\{\begin{array}{cl}
k&\textup{ if }n=2k,\\
0&\textup{ otherwise,}
\end{array}\right.\quad\textup{ and }\quad Sq^{\frac{n}{2}}=\left\{\begin{array}{cl}
Sq^{k}&\textup{ if }n=2k,\\
0&\textup{ otherwise}.
\end{array}\right.
$$
\begin{prop}[\cite{Sch94}]
	For all integers $n\geq 1,$ there are isomorphisms of unstable modules
	\begin{align*}
		\tilde{\Sigma}J(n)&\cong J(n-1),\\
		\tilde{\Phi}J(n)&\cong J\left(\frac{n}{2}\right).
	\end{align*}
	Moreover $\tilde{\lambda}_{J(n)}= \bullet Sq^{\frac{n}{2}}$.
\end{prop}
Recall that $\h^{*}\moc{B\z/2;\f2}$ is isomorphic to the polynomial algebra $\f2\left[u\right]$ on one variable of degree $1$, and the action of the Steenrod algebra is given by $Sq^{n}u^{k}=\binom{k}{n}u^{n+k}\moc{\mathrm{mod}\textbf{ }2}$. Since $F(1)$ can be identified as the submodule of $\h^{*}\moc{B\z/2;\f2}$ generated by $u$, then it has an $\f2-$basis consisting of $u^{2^{n}},n\geq 0.$ As a result, we have 
$$\ho{\u}{F(1)}{J\left(2^{n}\right)}\cong\f2$$
for all integers $n\geq 0$. Denote by $x_{n}$ the unique  generator in degree $1$ of $J\moc{2^{n}}$, and by $\mathscr{M}$ the bigraded algebra
$$\f2[x_{n},n\geq 0,||x_{n}||=(1,2^{n})].$$
The following theorem, due to Miller, describes the $J(n)$.
\begin{thm}[{\cite{Mil84}}]
	The morphism 
	\begin{align*}
		\bigoplus_{n\geq 0}J(n)\to\mathscr{M},\quad 
		x_{n}\mapsto x_{n}. 
	\end{align*}
	is an $\a2-$isomorphism of bigraded algebras.
\end{thm}
\begin{cor}
	Let $\sum\limits_{i=1}^{d}2^{j_{i}}$ be the $2-$adic expression of $n.$ Then, we have
	\begin{align*}
		J(n)^{n}&=\f2\left\langle x_{0}^{n}\right\rangle,\\
		J(n)^{d}&=\f2\left\langle \prod_{i=1}^{d}x_{j_{i}}\right\rangle,\\
		J(n)^{m}&=0\text{ if either } m>n \text{ or }d>m.
	\end{align*}
\end{cor}

\noindent \textsc{\bfseries The Bockstein long exact sequence. }Because the square of the Bockstein operation is trivial, the sequence $\left(J(n),\bullet
Sq^{1}:J(n+1)\to J(n)\right)^{n\geq 1}$ is a complex, which is close to being exact. In fact, we can modify this complex to obtain an exact sequence.
\begin{defi}[The Bockstein sequence]\label{bocksteins}
Let $\moc{\mathscr{B}_{k},\beta_{k}}$ be the following complex:
\begin{align*}
	\mathscr{B}_{-n}&=\left\{\begin{array}{cl}
	J(4k-1)\oplus J\left( 2k\right)&\textup{ if } n=4k-1,\textup{ and } n> 0,\\
	J(n)&\textup{ if } n\not\equiv 0(4),\textup{ and } n> 0,\\
	0&\textup{ otherwise.}
	\end{array} \right.\\ 
	\beta_{-n}&=\left\{\begin{array}{cl}
	\moc{\bullet Sq^{1}\quad 0}
	&\textup{ if } n=4k-1,\textup{ and } n> 0,\\
	\binom{\bullet Sq^{1}}{0}&\textup{ if } n=4k,\textup{ and } n> 0,\\
	0&\textup{ if } 0\geq n\\
	\bullet Sq^{1}&\textup{ otherwise.}
	\end{array} \right.
	\end{align*}
We call this complex \textit{the Bockstein sequence of Brown-Gitler modules}.
\end{defi}
\begin{prop}\label{bsss}
	The Bockstein sequence is exact.
\end{prop}
\begin{proof}
    We observe that an element $x\in J(n)$ belongs to $\ke{\bullet Sq^{1}:J(n)\to J(n-1)}$ if and only if there exists an admissible Steenrod operation
    	\begin{align*}
    		\theta&=Sq^{2k_{1}}Sq^{k_{2}}\ldots Sq^{k_{m}}
    	\end{align*}
    	such that $e(\theta)\leq |x|$, $ 4k_{1}<n,$ and $\theta x=x_{0}^{n}.$ Remark that $e(\theta)=|x|$ if and only if $n=4k_{1}$. Therefore,
    	\ali{e(Sq^{1}\theta)&\leq |x|\textup{ if }4\nmid n,\\
    	e(Sq^{1}\theta)&>|x|\textup{ otherwise.}
    	} 
    It means that if $n$ is not divisible by $4$, then 
    \ali{\ke{\beta_{n}}=\im{\beta_{n-1}}.}	
    It remains to show that 
    \ali{\ke{\beta_{-4k+1}}&=\im{\beta_{-4k}}}
    for all integers $k\geq 1$. This amounts to saying that all element of $J(2k)\subset \mathscr{B}_{4k-1}$ is a co-boundary. But this comes from the surjectivity of the map $\bullet Sq^{2k}:J(4k)\to J(2k)$, which is a direct consequence of Theorem \ref{mah}.
\end{proof}

\section{Graph representation of complexes}\label{graphrep}
In this section, we introduce the notion of \textit{graph representation of complexes}. In other words, we show how to associate a complex with an appropriate graph with respect to a certain decomposition. 
\subsection{Graph representation}
In all abelian category, a morphism between two direct sums of objects can be represented by matrix: let \ali{f:\dsum{i=1}{n}{M_i}\to \dsum{j=1}{m}{N_j}} 
be such a morphism, then we represent $f$ as an $m-$by$-n$ 
matrix 
$M_f=\left\{f_{ij}\right\},$
where $f_{ij}$ denotes the induced morphism $M_{j}\to N_{i}.$ It is sometimes convenient to respresent such a morphism as a bipartite graph with two disjoint sets of vertices $M,N$ indexed by $M_{i}$ and $N_{j}$ respectively, and the edges are determined by the matrix $M_{f}$.  
\begin{exam}
The morphism
\ali{\begin{pmatrix}
\bullet Sq^{1}& 0\\
\bullet Sq^{2,1}& \bullet Sq^{2}\\
0 & \bullet Sq^{3} 
\end{pmatrix}:J(7)\oplus J(6)\longrightarrow J(6)\oplus J(4)\oplus J(3)}
is represented as follows:
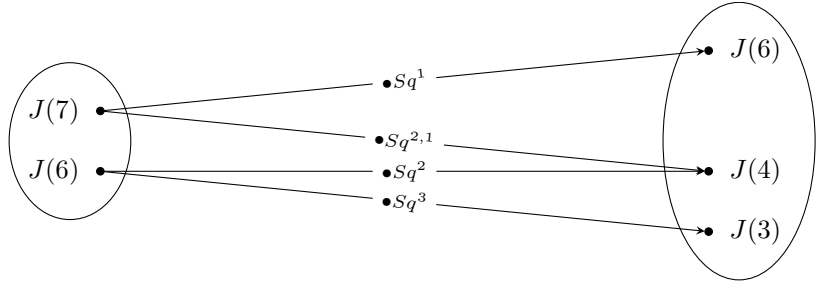
\begin{center}
\input{exam}
\end{center}
\end{exam}
\begin{defi}
Let $\moc{C^{k},\partial^{k}}$ be a cochain complex. Suppose that each $C^{k}$ admits a decomposition
\alilabel{decomposition}{C^{k}\cong \bigoplus_{\alpha\in A_{k}}C^{k}_{\alpha}}
and denote by $M_{\partial^{k}}$ the corresponding matrix representing $\partial^{k}$. Then, the graph associated with $\moc{C^{k},\partial^{k}}$ with respect to the decomposition \eqref{decomposition} is defined as follows:
\begin{enumerate}
\item The set of vertices is the disjoint union of $V_{k}$ indexed by $A_{k}$:
\ali{V_{k}:=\left\{\left.  v_{\alpha}^{k} \right| \alpha\in A_{k}  \right\}.}
\item The set of edges is the disjoint union of $E_{k}$ containing edges from $V_{k}$ to $V_{k+1}$ determined by the matrix $M_{\partial^{k}}$:
\ali{E_{k}:=\left\{\left. \left[ v_{\alpha}^{k},v_{\beta}^{k+1}  \right]\right| 0\neq\partial^{k}_{\beta,\alpha}:C_{\alpha}^{k}\to C_{\beta}^{k+1}     \right\}.}
\end{enumerate}
We write $\left[ v_{\alpha}^{k},v_{\beta}^{k+1}  \right]=\partial^{k}_{\beta,\alpha}$, with the convention that $\left[ v_{\alpha}^{k},v_{\beta}^{k+1}  \right]=0$ means there is no edge between  $v_{\alpha}^{k}$ and $v_{\beta}^{k+1}$.  
\end{defi}
\subsection{Graph representation of resolutions}
We now consider some particular complexes of unstable modules and their graph representation.

\begin{defi}[BG modules and BG complexes]
A BG module is a direct sum of Brown-Gitler modules, and a complex of BG modules is called a BG complex. 
\end{defi}
\begin{defi}
Let $f:\dsum{\alpha\in A}{}{J\moc{n_{\alpha}}}\to \dsum{\beta\in B}{}{J\moc{m_{\beta}}}$ be a morphism between BG modules, and let 
$$M_f=\left\{f_{\beta,\alpha}:J\moc{n_{\alpha}}\to J\moc{m_{\beta}}\right\}$$
be its representing matrix. We define $\bar{f}$ to be the morphism with the same source and target as those of $f$, but the representing matrix is obtained from that of $f$ by identifying all $f_{\beta,\alpha}$ with $0$ whenever $f_{\beta,\alpha}$ is not the identity map $\bullet Sq^{0}$.
\end{defi}
The following lemma is straightforward.
\begin{lem}\label{sq0}
Let $\moc{I^{\bullet},\partial^{\bullet}}$ be a BG complex, then so is $\moc{I^{\bullet},\bar{\partial^{\bullet}}}$.
\end{lem}
\begin{defi}\label{graphdeg}
Under the same hypothesis as that of Lemma \ref{sq0}, we denote by $G\moc{I^{\bullet},\bar{\partial^{\bullet}}}$ the associated graph with respect to the decomposition of $I^{\bullet}$ as direct sum of Brown-Gitler modules. We define  $R\moc{I^{\bullet},\partial^{\bullet}}$ to be the bigraded $\f2-$vector space generated by the vertices of $G\moc{I^{\bullet},\bar{\partial^{\bullet}}}$, where $v\in V_{r}\moc{G\moc{I^{\bullet},\bar{\partial^{\bullet}}}}$ is of bidegree $(r,s)$ if $v$ is indexed by a direct summand of the form $J(s)$. We call this the \textit{BG degree} of $v$ and denote it by $\left|v\right|=\moc{\left|v\right|_{1},\left|v\right|_{2}}$, where $\left|v\right|_{1}=r$ and $\left|v\right|_{2}=s$.
\end{defi}
We now show how to relate Ext-groups of unstable modules to the graph representation of resolutions. The following lemma is the key to our construction of the Lambda algebra in Section \ref{Lambdaalgebra}.
\begin{lem}\label{lambdaext}
Let $M$ be a finite unstable module. Then, $M$ admits an injective resolution $\moc{I^{\bullet},\partial^{\bullet}}$ of finite length such that each $I^{k}$ is a finite BG module for all $k\geq 0$. Moreover, we have
\alilabel{extgroup1}{\h^{r,s}\moc{R\moc{I^{\bullet},\partial^{\bullet}}}\cong \ext{\u}{r}{\Sigma^{s}\f2}{M}}
for all integers $r,s\geq 0.$
\end{lem}
\begin{proof}
Denote by $d\moc{M}$ the least degree $n$ such that $M^{k}=0$ for all $k>n$. The existence of such an injective resolution $\moc{I^{\bullet},\partial^{\bullet}}$ can be proved by induction on the degree $d\moc{M}$. We now show that the isomorphism \eqref{extgroup1} holds. Indeed, for all integers $r,s\geq 0$, we have
\ali{\ext{\u}{r}{\Sigma^{s}\f2}{M}&\cong \h^{r}\ho{\u}{\Sigma^{s}\f2}{\moc{I^{\bullet},\partial^{\bullet}}}\\
&\cong \h^{r}\ho{\u}{\Sigma^{s}\f2}{\moc{I^{\bullet},\bar{\partial^{\bullet}}}}\\
&\cong \h^{r,s}\moc{R\moc{I^{\bullet},\partial^{\bullet}}},}
whence the conclusion.
\end{proof}

\section{Pseudo-hyperresolutions}\label{Phr}
Constructing resolutions is one of the most basic problems in homological
algebra. This section aims to study a certain class of objects whose resolutions can be made explicit. We begin with some elementary examples which are the origin of the present paper. 

Recall that the groups $\exg{*}{m}{n}$ are of interest because of the unstable Adams spectral sequence computing the $2-$component of the homotopy groups of spehres. Therefore, it is natural to search for explicit injective resolutions of $\Sigma^{n}\f2$ in the category $\u$. The cases $n=0$ and $n=1$ are trivial since $\Sigma^{n}\f2\cong J(n)$ in these cases. As a result, the Mahowald short exact sequence
\alilabel{sigma2}{0\to\Sigma J(1)\xrightarrow{s_{1}} J(2)\xrightarrow{\bullet Sq^{1}}J(1)\to 0}
is an injective resolution of $\Sigma^{2}\f2$. In order to construct an injective resolution for $\Sigma^{3}\f2$, we apply the suspension functor $\Sigma$ to the sequence \eqref{sigma2} and then obtain another exact sequence. But, this is no longer a resolution. However, we have the following commutative diagram.
\alilabel{doublecomplex}{\input{sigma2}} 
Consequently, the total complex of the commutative square \eqref{doublecomplex} is an injective resolution of $\Sigma^{3}\f2$. 

The above examples lead to the following question.

\begin{quest}\label{quest1}
Let $M$ be an object in an abelian category $\mathcal{C}$ such that:
\begin{itemize}
	\item The category $\mathcal{C}$ has enough projective (injective) objects.
	\item There exists an acyclic complex in $\mathcal{C}$ such that its only nontrivial (co-)homology is
	isomorphic to $M,$ and that each term of the complex admits an explicit projective (injective)
	resolution. 
\end{itemize}
Is it possible to construct an explicit resolution for $M$?
\end{quest}

The construction of an explicit resolution for such an object $M$ is the main goal of this section, and the results are stated in Propositions \ref{key121} and \ref{key12}. 

Throughout this section, a complex refers to a cohomological one unless otherwise stated. 
\subsection{The local-to-global principle}
The construction of an explicit resolution for such an object $M$ in Question \ref{quest1} will be carried out by a \textit{local-to-global} principle. In other words, for all integers $n\geq 0$, we  construct an acyclic complex $\moc{I^{k}(n),\partial^{k}(n)}$ such that:
\begin{enumerate}
\item The only nontrivial cohomology of this complex is $M$.
\item The term $I^{k}(n)$ is injective for all $k\leq n$.
\item We have $I^{k}(n)=I^{k}(n+1)$ for all $k\leq n$.
\end{enumerate}
Therefore, by letting $n$ tend to infinity, we obtain an injective resolution of $M$. Before formulating the principle, we need the following key lemma.
\begin{lem}\label{bg2}
	In an abelian category, if the sequence
	$$E_1\xrightarrow{f_1}E_2\oplus M\xrightarrow{\left( \begin{smallmatrix}
		g_1 & g_2\\
		g_3 & id\end{smallmatrix} \right)}E_3\oplus M\xrightarrow{f_3}E_4$$
	is exact at $E_2\oplus M$ and $E_3\oplus M$, then the sequence 
	$$E_1\xrightarrow{f_1}E_2\xrightarrow{g_1+g_2\circ g_3}E_3\xrightarrow{f_3}E_4$$
	is exact at $E_2$ and $E_3$. 
\end{lem}
\begin{proof}
	Note that the diagram 
	$$\xymatrixcolsep{4pc}\xymatrixrowsep{2pc}\xymatrix{0
		\ar[r]^{}\ar[d]&M\ar[r]^{id}\ar[d]^{\binom{0}{id}}&M\ar[r]^{}\ar[d]^{\binom{g_2}{id}}& 0\ar[d]\\
		E_1\ar[r]^{f_1}\ar[d]^{id}& E_2\oplus M\ar[r]^{\left( \begin{smallmatrix}
			g_1 & g_2\\
			g_3 & id\end{smallmatrix} \right)}\ar[d]^{(id,0)}&E_3\oplus
		M\ar[r]^{f_3}\ar[d]^{(id,g_2)}&E_4\ar[d]^{id}\\
		E_1\ar[r]^{f_1}&E_2\ar[r]^{g_1+g_2\circ g_3}&E_3\ar[r]^{f_3}&E_4
	}$$
	is commutative, and each of its columns is exact. As the first two upper rows are exact at the second and the third terms, then so is the last row.
\end{proof}
The argument now goes as follows.
\begin{lem}\label{key9}
	In an abelian category with enough injective objects, let
	$\left(M^{k},\tau^{k}:M^{k}\to M^{k+1},k\geq 0\right)$ be an acyclic complex, and let $M$ be
	its only nontrivial cohomology. Denote by $I^{k}$ an injective object containing $M^{k}$ and by
	$p^{k}$ the induced projection $I^{k}\to I^{k}/M^{k}.$ Then, for all integers $k\geq 0,$ there exist morphisms
	\ali{
		\theta^{k}: I^{k}\to I^{k+1},\quad
		\delta^{k}: \frac{I^{k}}{M^{k}}\to I^{k+2},\quad
		\omega^{k}: \frac{I^{k}}{M^{k}}\to \frac{I^{k+1}}{M^{k+1}},}
	such that the morphisms $\left(\begin{smallmatrix}
		\theta^k &\delta^{k-1}\\
		p^{k} & \omega^{k-1}
		\end{smallmatrix}\right)$ for all integers $k\geq 1$, denoted by
	$\partial^{k}$, and $\left(\begin{smallmatrix}
		\theta^0\\
		p^{0}
		\end{smallmatrix}\right)$, denoted by $\partial^{0}$, 
	make the sequence 
	\begin{equation}\label{suite1}
		0\to M\to I^{0}\xrightarrow{\partial^{0}} I^{1}\bigoplus
		\frac{I^{0}}{M^{0}}\xrightarrow{\partial^{1}}\cdots\xrightarrow{\partial^{t}}I^{t+1}\bigoplus
		\frac{I^{t}}{M^{t}}\xrightarrow{\partial^{t+1}}\cdots
	\end{equation}
	an exact complex.
\end{lem}
\begin{proof}
	For all integers $k\geq 0,$ denote by $i^{k}$ the inclusion $M^{k}\to I^{k}$. Because $I^{k}$
	is injective for all integers $k\geq 0$, the $\tau^{k}$ gives rise to a morphism $\theta^{k}:I^{k}\to I^{k+1}$ such
	that 
	$$\theta^{k}\circ
	i^{k}=i^{k+1}\tau^{k}.$$
	As $\theta^{k+1}\circ\theta^{k}\circ i^{k}$ is trivial for all integers $k\geq 0$, there exists a morphism
	$\delta^{k}: I^{k}/M^{k}\to I^{k+2}$ such that 
	$$\delta^{k}\circ p^{k}=\theta^{k+1}\circ\theta^{k}.$$
	Denote by 
	$$\omega^{k}:\frac{I^{k}}{M^{k}}\to \frac{I^{k+1}}{M^{k+1}}$$ 
	the morphism induced by $\theta^{k}$ for all integers $k\geq 0$.
	Then, $\left(p^{k},k\geq 0\right)$ is a morphism of complexes: 
	$$\left(I^{k},\theta^{k},k\geq
	0\right)\to \left(\frac{I^{k}}{M^{k}},\omega^{k},k\geq 0\right).$$
	Denote by $K^{l}$ the direct sum $I^{l}\oplus I^{l+1}$ for all integers $l\geq 1$, and by $K^{0}$ the module $I^{0}$.
	Then, the morphisms
	\begin{align*}
		\beta^0&=\left( \begin{smallmatrix}
			\theta^0\\
			\theta^1\circ\theta^0\end{smallmatrix} \right),\\ 
		\beta^l&=\left( \begin{smallmatrix}
			\theta^l & id\\
			\theta^{l+1}\circ\theta^l & \theta^{l+1}\end{smallmatrix} \right),\quad l\geq 1,
	\end{align*}
	make $\left(K^{l},\beta^{l}:K^{l}\to K^{l+1},l\geq 0\right)$ a complex.
		Denote by $H^{l}$ the direct sum $I^{l+1}\oplus I^{l}/M^{l}$ for all integers $l\geq 1$, and denote by $H^{0}$ the module
	$I^{0}/M^{0}$.  Then, the morphisms
	\begin{align*}
		\gamma^0&=\left( \begin{smallmatrix}
			\theta^0\\
			\omega^0\end{smallmatrix} \right),\\
		\gamma^k&=\left( \begin{smallmatrix}
			\theta^{k+1} & \delta^k\\
			p^{k+1} & \omega^k\end{smallmatrix} \right),\quad k\geq 1,
	\end{align*}
	make $\left(H^{k},\gamma^{k},k\geq 0\right)$ a complex. For all integers $k\geq 1$, we fix:
	$$
	\begin{aligned}
	\eta^0&=i^0,& \eta^k&=\left( \begin{smallmatrix}
	i^{k}\\
	0\end{smallmatrix} \right),&\\
	\chi^0&=p^0,& \chi^k&=\left( \begin{smallmatrix}
	0 & id\\
	p^{k} & 0\end{smallmatrix} \right).&
	\end{aligned}$$
	Then, the sequence
	$$0\to M^{k}\xrightarrow{\eta^{k}} K^{k}\xrightarrow{\chi^{k}}H^{k}\to 0$$
	is exact for all integers $k\geq0.$ Since $\left(\eta^{l},l\geq 0\right)$ and  $\left(\chi^{l},l\geq 0\right)$ are morphisms of complexes we obtain the double complex
	\alilabel{doublecomplex}{\xymatrix{
		K^{0}\ar[r]^{\beta^0}\ar[d]^{\chi^0}&
		K^{1}\ar[r]^{\beta^1}\ar[d]^{\chi^1}&\cdots\ar[r]^{\beta^{n-2}}& K^{n-1}\ar[r]^{
			\beta^{n-1}}\ar[d]^{\chi^{n-1}}&\cdots\\
		H^{0}\ar[r]^{\gamma^0}& H^{1}\ar[r]^{\gamma^1}&\cdots\ar[r]^{\gamma^{n-2}} & H^{n-1}
		\ar[r]^{
			\gamma^{n-1}} & \cdots}}
	where all the columns are acyclic. As the sequence $\left(M^{k},\tau^{k}:M^{k}\to M^{k+1},k\geq 0\right)$ is also acyclic, using standard spectral sequence arguments for double complexes, it is straightforward that the total complex of the double complex \eqref{doublecomplex} is acyclic and its only nontrivial cohomology is $M$. By applying Lemma \ref{bg2} to this
	acyclic sequence, we obtain the exact sequence \eqref{suite1}.
\end{proof}
By the same method, we can prove the following generalized version of Lemma \ref{key9}, which does not require complexes to be bounded below.
\begin{cor}\label{corkey1}
		In an abelian category with enough injectives, let
		$\left(M^{k},\tau^{k}:M^{k}\to M^{k+1}\right)$ be a complex. Denote by $I^{k}$ an injective object containing $M^{k}$, and denote by
		$p^{k}$ the induced projection $I^{k}\to I^{k}/M^{k}.$ Then, for all integers $k$, there exist morphisms
		\begin{align*}
		\theta^{k}: I^{k}\to I^{k+1},\quad
		\delta^{k}: \frac{I^{k}}{M^{k}}\to I^{k+2},\quad
		\omega^{k}: \frac{I^{k}}{M^{k}}\to \frac{I^{k+1}}{M^{k+1}},
		\end{align*}
		such that the morphisms $\left(\begin{smallmatrix}
				\theta^{k} &\delta^{k-1}\\
				p^{k} & \omega^{k-1}
				\end{smallmatrix}\right)$, denoted by 
		$\partial^{k}$, 
		make the sequence 
		\begin{equation}
		\cdots\xrightarrow{\partial^{t-1}} I^{t}\bigoplus \frac{I^{t-1}}{M^{t-1}}\xrightarrow{\partial^{t}} I^{t+1}\bigoplus
		\frac{I^{t}}{M^{t}}\xrightarrow{\partial^{t+1}}\cdots
		\end{equation}
		a complex. Moreover, the compositions \ali{M^{k}\hookrightarrow I^{k}\xrightarrow{\binom{id}{0}} I^{k}\bigoplus I^{k-1}/M^{k-1}} 
		form a morphism of complexes, which induces isomorphisms on cohomology groups.
\end{cor}
\subsection{The construction}
The following proposition is a strengthening of Lemma \ref{key9}.
\begin{prop}\label{key121}
	In an abelian category with enough injectives, let
	$\left(M^{l},\tau^{l}:M^{l}\to M^{l+1},l\geq 0\right)$ be an acyclic complex, and denote by $M$
	its only nontrivial cohomology. For all $l\geq 0$, let $\left(I^{l,t},\partial^{l,t}:I^{l,t}\to I^{l,t+1},k\geq t\geq 0\right)_{}$ be an acyclic complex such that its only nontrivial cohomology is $M^{l},$ and $I^{l,t}$
	is injective for all $k> t\geq 0$, and $I^{l,k}$ is the quotient
	$I^{l,k-1}/\im{\partial^{l,k-2}}$. Then, there exist morphisms
	$$\partial^{l,t}_{i}:I^{l,t}\to
	I^{l+i+1,t-i}$$ 
	such that the morphisms
	$$\partial^{l}=\begin{pmatrix}
	\partial^{0,l} & 0 & \cdots & 0 \\
	\partial^{0,l}_{0} & \partial^{1,l-1} & \cdots & 0 \\
	\vdots  & \vdots  & \ddots & \vdots  \\
	\partial^{0,l}_{l-1} & \partial^{1,l-1}_{l-2} & \cdots & \partial^{l,0}\\
	\partial^{0,l}_{l} & \partial^{1,l-2}_{l-1} & \cdots & \partial^{l,0}_{0}
	\end{pmatrix}$$
	make the sequence
	\begin{equation}\label{eqk}
	0\to M\to  I^{0,0}\xrightarrow{\partial^{0}} I^{0,1}\bigoplus
	I^{1,0}\xrightarrow{\partial^{1}}\cdots\xrightarrow{\partial^{l}}\dsum{\substack{m+n=l+1\\0\leq
				n\leq k}}{}{I^{m,n}}
	\xrightarrow{\partial^{l+1}}\cdots
	\end{equation}
	exact.
\end{prop}
\begin{proof}
	We proceed by induction on $k$. The case $k=1$ is proved in Lemma \ref{key9}. Suppose that the
	lemma is true for $k<q$. We show that it is also true for $k=q$. Denote by
	$J^{i,q-1}$ the quotient $I^{i,q-2}/\im{\partial^{i,q-3}}.$ By induction hypothesis, there exists
	an exact sequence:
	$$M\hookrightarrow I^{0,0}\xrightarrow{} I^{0,1}\bigoplus
	I^{1,0}\xrightarrow{}\cdots\xrightarrow{}\left(\dsum{\substack{m+n=l+1\\0\leq
			n<q-1}}{}{I^{m,n}}\right)\oplus J^{l-q+2,q-1}\xrightarrow{}\cdots$$
	Because the sequence 
	$$0\to\left(\dsum{\substack{m+n=l+1\\0\leq n<q-1}}{}{I^{m,n}}\right)\oplus
	J^{l-q+2,q-1}\to \left(\dsum{\substack{m+n=l+1\\0\leq
			n<q-1}}{}{I^{m,n}}\right)\oplus
	I^{l-q+2,q-1}\to I^{l-q+2,q}\to 0$$
	is exact, then, according to Lemma \ref{key9}, the complex
	\begin{equation*}
		0\to M\to I^{0,0}\xrightarrow{\partial^{0}} I^{0,1}\bigoplus
		I^{1,0}\xrightarrow{\partial^{1}}\cdots\xrightarrow{\partial^{l}}\dsum{\substack{m+n=l+1\\0\leq
						n\leq k}}{}{I^{m,n}}
		\xrightarrow{\partial^{l+1}}\cdots
	\end{equation*}
	is exact. We can then conclude the proposition.
\end{proof}
Note that, in the exact sequence \eqref{eqk}, the first $k-1$ terms are injective. Let $k$ tend to infinity, we
obtain an injective resolution of $M$.
\begin{depo}[Pseudo-hyperresolution]\label{key12}
	In an abelian category with enough injective objects, let
		$\left(M^{k},\tau^{k}:M^{k}\to M^{k+1},k\geq 0\right)$ be an acyclic complex, and let $M$ be
		its only nontrivial cohomology. For all integers $k\geq 0$, let $\left(I^{k,t},\partial^{k,t}\right)_{t\geq
		0}$ be an injective resolution of $M^{k}.$ Then, there exist morphisms
	$$\partial^{k,t}_{i}:I^{k,t}\to
	I^{k+i+1,t-i}$$ such that the morphisms
	$$\partial_{k}=\begin{pmatrix}
	\partial^{0,k} & 0 & \cdots & 0 \\
	\partial^{0,k}_{0} & \partial^{1,k-1} & \cdots & 0 \\
	\vdots  & \vdots  & \ddots & \vdots  \\
	\partial^{0,k}_{k-1} & \partial^{1,k-1}_{k-1} & \cdots & \partial^{k,0}\\
	\partial^{0,k}_{k} & \partial^{1,k-2}_{k} & \cdots & \partial^{k,0}_{0}
	\end{pmatrix}$$
	make the complex
	$$I^{0,0}\xrightarrow{\partial^{0}} I^{0,1}\bigoplus
	I^{1,0}\xrightarrow{\partial^{1}}\cdots\xrightarrow{\partial^{k}}\dsum{m+n=k+1}{}{I^{m,n}}
	\xrightarrow{\partial^{k+1}}\cdots$$
	an injective resolution of $M$. This resolution, denoted by $\psr{M}{I^{\bullet,\bullet}}$, is called the \phr of $M$ with respect to the resolutions 
	$\left(I^{k,t},\partial^{k,t}\right)_{t\geq
			0}$.
\end{depo}
Because of Corollary \ref{corkey1}, we also have a generalized version of Proposition \ref{key12}.
\begin{cor}\label{corkey2}
		In an abelian category with enough injectives, let
				$\left(M^{k},\tau^{k}:M^{k}\to M^{k+1}\right)$ be a complex.  Moreover, for all integers $k$, let $\left(I^{k,t},\partial^{k,t}\right)_{t\geq 0}$ be an injective resolution of $M^{k}$. Then, there exist morphisms
		$$\partial^{k,t}_{i}:I^{k,t}\to
		I^{k+i+1,t-i}$$ such that the morphisms
		$$\partial_{k}=\begin{pmatrix}
		\partial^{0,k} & 0 & \cdots & 0 \\
		\partial^{0,k}_{0} & \partial^{1,k-1} & \cdots & 0 \\
		\vdots  & \vdots  & \ddots & \vdots  \\
		\partial^{0,k}_{k-1} & \partial^{1,k-1}_{k-2} & \cdots & \partial^{k,0}\\
		\partial^{0,k}_{k} & \partial^{1,k-2}_{k-1} & \cdots & \partial^{k,0}_{0}
		\end{pmatrix}$$
		make the following sequence
		$$ \cdots\xrightarrow{\partial^{k-1}}\dsum{m+n=k}{}{I^{m,n}}\xrightarrow{\partial^{k}}\dsum{m+n=k+1}{}{I^{m,n}}
		\xrightarrow{\partial^{k+1}}\cdots$$
		a complex, and the compositions \ali{M^{k}\hookrightarrow I^{k,0}\hookrightarrow \dsum{m+n=k}{}{I^{m,n}}}
		form a morphism of complexes, which induces isomorphisms on cohomology groups.
\end{cor}
Therefore, in all abelian categories, each cochain complex is quasi-isomorphic to a complex of injective objects. Here, recall that two complex are quasi-isomorphic if there exists a morphism between them that induces isomorphisms on cohomology groups.
\begin{rem}
The technical point we used in the proof of Lemma \ref{key9} is similar to the mapping telescope technique in triangulated categories introduced by B\"{o}ckstedt and Neeman \cite{BN93}. Then, the intrigued readers might ask what is the advantage of introducing the notion of a pseudo-hyperresolution? The answer lies in the use of Lemma \ref{bg2}: this makes our method different from that of B\"{o}ckstedt and Neeman. In fact, when the length of the complex in play is infinite, the resulting resolution we obtained is more compact and, hence, more computable. To illustrate the value of the pseudo-hyperresolution method, we will give an application for infinite-length complexes in Section \ref{icps}.
\end{rem}

Dually, we can construct the so-called pseudo-hyper projective resolutions.
\begin{cor}\label{phpro}
In an abelian category with enough projective objects, let
	$\left(N^{k},\tau^{k+1}:N^{k+1}\to N^{k},k\geq 0\right)$ be an acyclic complex, and denote by $N$ be
	its only nontrivial homology. For all $k\geq 0$, let $\left(P_{k,t},\partial_{k,t+1}:P_{k,t+1}\to P_{k,t}\right)_{t\geq
		0}$ be a projective resolution of $N^{k}.$ There exist morphisms  $$\partial_{k,t}^{i}:P_{k,t}\to
	P_{k-i-1,t+i}$$ 
	such that the morphisms
	$$\partial_{k}=\begin{pmatrix}
	\partial_{0,k} & \partial_{1,k-1}^{0} & \cdots & \partial_{k-1,1}^{k-2} & \partial_{k,0}^{k-1}\\
	0 & \partial_{1,k-1} & \cdots & \partial_{k-1,1}^{k-3} & \partial_{k,0}^{k-2}\\
	\vdots  & \vdots  & \ddots & \vdots &\vdots \\
	0 & 0 & \cdots & \partial_{k-1,1} & \partial_{k,0}^{0}
	\end{pmatrix}$$
	make the following sequence
	$$\cdots\xrightarrow{\partial_{k+2}}\dsum{m+n=k+1}{}{P_{m,n}}\xrightarrow{\partial_{k+1}}\cdots
	\xrightarrow{\partial_{2}}P_{0,1}\bigoplus P_{1,0} \xrightarrow{\partial_{1}} P_{0,0}$$
	a projective resolution of $N$. This resolution, denoted by $\psr{N}{P_{\bullet,\bullet}}$, is called the \phr of $N$ with respect to the resolutions 
	$\left(P_{k,t},\partial_{k,t+1}:P_{k,t+1}\to P_{k,t}\right)_{t\geq
			0}$.
\end{cor}

Remark that each term of the \phr of an object in an abelian category is explicit. Howerver, in general, it is difficult to determine the differentials of this resolution. Therefore, it is convenient to the notion of a \textit{Poincar\'{e} series} of a complex. 
\begin{defi}[Poincar\'{e} series of complexes]
In an abelian category, let $\moc{C_{\bullet},\partial_{\bullet}}$ be a complex. Then, we denote by $\h\moc{t,C_{\bullet}}$ the series
\alilabel{poincare}{\h\moc{t,C_{\bullet}}:=\sum_{i}^{}\mocv{C_{i}}t^{i},}
where $\mocv{C_{i}}$ denotes the corresponding class of $C_{i}$ in the Grothendieck group of the category. 
\end{defi} 
The following lemmas are straightforward.
\begin{lem}
In an abelian category, let $\moc{C_{\bullet},\partial_{\bullet}}$ and $\moc{D_{\bullet},\delta_{\bullet}}$ be two complexes. Then, we have
\ali{\h\moc{t,C_{\bullet}\oplus D_{\bullet}}&= \h\moc{t,C_{\bullet}}+\h\moc{t,D_{\bullet}}.}
Moreover, if the category is monoidal, then we have
\ali{
\h\moc{t,C_{\bullet}\otimes D_{\bullet}}&= \h\moc{t,C_{\bullet}}\cdot\h\moc{t,D_{\bullet}}.}
\end{lem}
\begin{lem}
Under the same hypothesis as that of Corollary \ref{phpro}, we have
\ali{
\h\moc{t,\psr{N}{P_{\bullet,\bullet}}}&= \sum_{r,s\geq 0}^{}P_{r,s}t^{r+s}.}
\end{lem}

\section{Functor homology}\label{functor}

The main purpose of this section is to provide two simple applications of the \phr method in functor homology:
\begin{enumerate}
\item Constructing explicit resolutions of certain strict polynomial functors and use them to compute the Mac Lane cohomology of $\f2$ (see Theorem \ref{mclane}).
\item Computing the global dimension of the category of homogeneous strict polynomial functors of finite degree (see Theorem \ref{globaldimension}). 
\end{enumerate}

\subsection{Mac Lane cohomology of finite fields}
Following the works \cite{EM53,EM54,EM54b} of Eilenberg and Mac Lane on the homology of the spaces that bear their names, Mac Lane introduces in \cite{Mac57} the notion of the  cohomology of a ring with coefficients in a bimodule. Jibladze and Pirashvili extend this notion for more general coefficients (see \cite{JP91}). In what follow, we recall the definition of Mac Lane cohomology and show how to use the \phr method to compute the Mac Lane cohomology of the finite field $\f2$.

\subsubsection{Ordinary functors and strict polynomial functors}
Let $R$ be a ring and denote by $\mathscr{F}(R)$ the category of functors from the category of projective $R-$modules of finite rank to the category of $R-$modules. Denote by $I$ the inclusion functor that associates a  projective $R-$modules of finite rank $V$ with itself. Then, following \cite{JP91}, we define the Mac Lane cohomology of $R$ with coefficients in a functor $F\in\mathscr{F}(R)$ as follows:
\ali{\hml{R}{F}:=\ext{\mathscr{F}(R)}{*}{I}{F}.}
Taking $F=I$, we recover the original  Mac Lane cohomology of the ring $R$ with coefficients in $R$, with the obvious structure of $R-R-$bimodule (see \cite{Mac57}).

In fact, Ext-groups in $\mathscr{F}(R)$ can capture more than just the Mac Lane cohomology. When $R$ is a finite field $\Bbbk$ of characteristic $p$, we write $\fk$ for $\mathscr{F}(\Bbbk)$. If $F\in\fk $, then structure morphism
\ali{F_{\Bbbk^{n},\Bbbk^{n}}:\ho{\Bbbk}{\Bbbk^{n}}{\Bbbk^{n}}\to\ho{\Bbbk}{F\moc{\Bbbk^{n}}}{F\moc{\Bbbk^{n}}}}
turns $F\moc{\Bbbk^{n}}$ into a $\Bbbk\left[\gln{n}{\Bbbk}\right]-$module. Therefore, the evaluation on $\Bbbk^{n}$ yields a natural morphism
\ali{\ext{\fk}{*}{F}{G}\to \ext{\Bbbk\left[\gln{n}{\Bbbk}\right]}{*}{F\moc{\Bbbk^{n}}}{G\moc{\Bbbk^{n}}}.}
The canonical inclusion 
\ali{\gln{n}{\Bbbk}\to \gln{n+1}{\Bbbk},\quad M\mapsto \begin{pmatrix}
M&0\\0&1
\end{pmatrix},}
along with the splitting projection $\Bbbk^{n+1}\to \Bbbk^{n}$ onto the first $n$ coordinates, induces a natural map
\ali{\ext{\Bbbk\left[\gln{n+1}{\Bbbk}\right]}{*}{F\moc{\Bbbk^{n+1}}}{G\moc{\Bbbk^{n+1}}}\to \ext{\Bbbk\left[\gln{n}{\Bbbk}\right]}{*}{F\moc{\Bbbk^{n}}}{G\moc{\Bbbk^{n}}}}
As $n$ tends to infinity, this map stabilizes, and we denote by $\ext{\Bbbk\left[\gln{\infty}{\Bbbk}\right]}{*}{F}{G}$ the stable value. Under some mild conditions on $F$ and $G$ (namely, $F,G$ are polynomial in the sense of Eilenberg and Mac Lane \cite{EM54}), the induced map
\ali{\ext{\fk}{*}{F}{G}\to \ext{\Bbbk\left[\gln{\infty}{\Bbbk}\right]}{*}{F}{G}}
is an isomorphism (see, \textit{e.g.}, \cite{Bet99,FFSS99}). Recall that if $M,N$ are $\Bbbk\left[\gln{n}{\Bbbk}\right]-$modules, then $\ho{\Bbbk}{M}{N}$ is endowed with an action of $\gln{n}{\Bbbk}$, and, moreover, there is an isomorphism
\ali{\ext{\Bbbk\left[\gln{\infty}{\Bbbk}\right]}{*}{F}{G}\cong \h^{*}\moc{\gln{n}{\Bbbk},\ho{\Bbbk}{M}{N}}.}
So, Ext-groups in $\fk$ can capture the stable cohomology of $\gln{n}{\Bbbk}$. Now, if we consider the affine algebraic group scheme $\Gln{n}{\Bbbk}$ instead, \textit{i.e.}, the functor that associates a $\Bbbk-$algebra $A$ with $\gln{n}{A}$, then it is natural to ask whether there exists a suitable category that is related to $\Gln{n}{\Bbbk}-$modules the way $\fk$ does to $\gln{n}{\Bbbk}-$modules. In \cite{FS97}, Friedlander and Suslin introduce the category $\p_{\Bbbk}$ of strict polynomial functors, resolving in the  affirmative this question. A strict polynomial functor $F$ is a functor from the category of $\Bbbk-$vector spaces of finite dimension to the category of $\Bbbk-$vector spaces such that the structure morphism
\alilabel{strucmap}{F_{V,W}:\ho{\Bbbk}{V}{W}\to\ho{\Bbbk}{F\moc{V}}{F\moc{W}}}
is a polynomial map. Here, a polynomial map between two $\Bbbk-$vector spaces $X$ and $Y$ is a morphism between two schemes $Spec\moc{S^{*}\moc{X^{\sharp}}}$ and $Spec\moc{S^{*}\moc{Y^{\sharp}}}$, where $S^{*}\moc{X^{\sharp}}$ is the symmetric algebra on the $\Bbbk-$linear dual  $X^{\sharp}$. Let $Sch/\Bbbk$ denote the category of schemes over $\Bbbk$. Because
\ali{\ho{Sch/\Bbbk}{Spec\moc{S^{*}\moc{X^{\sharp}}}}{Spec\moc{S^{*}\moc{Y^{\sharp}}}}\cong S^{*}\moc{X^{\sharp}}\otimes Y,}
then we define a polynomial map $p:X\to Y$ to be homogeneous of degree $d$ if it belongs to $S^{d}\moc{X^{\sharp}}\otimes Y$. And then, a strict polynomial functor is homogeneous of degree $d$ if its structure morphisms \eqref{strucmap} are homogeneous polynomial maps of the same degree. Denote by $\p_{\Bbbk,d}$ the full subcategory of homogeneous strict polynomial functors of degree $d$. Then,  $\p_{\Bbbk}\cong\bigoplus_{d\geq 0}\p_{\Bbbk,d}$.

As the evaluation on $\Bbbk$ of a morphism of schemes over $\Bbbk$ yields a set theoretic map, a strict polynomial functor $F$, with structural morphisms $\moc{F(V),F_{V,W}}$ gives rise to an ordinary functor $F\in\fk$, with structural morphisms $\moc{F(V),F_{V,W}(\Bbbk)}$. We obtain in this way an exact forgetful functor $\mathcal{O}:\p_{\Bbbk} \to\fk$.

Let $G\in Sch/\Bbbk$, then a representation of $G$ (or a $G-$module) is a $\Bbbk-$vector space, endowed with a natural transformation  $G\to\mathrm{GL}_{M}$, where $\mathrm{GL}_{M}$ is the functor that associates a $\Bbbk-$algebra $A$ with the group  $\mathrm{GL}_{A}\moc{A\otimes_{\Bbbk}M}$ of invertible $A-$linear endomorphisms of $A\otimes_{\Bbbk}M$. Evaluation on $\Bbbk$ yields an exact forgetful functor $\mathcal{O}:\Gln{n}{\Bbbk}-\mathrm{mod}\to\gln{n}{\Bbbk}-\mathrm{mod}$, and we have a commutative diagram
\alilabel{diagram3}{\xymatrix{
\p_{\Bbbk}\ar[d]_{\mathrm{ev}}\ar[r]^{\mathcal{O}}&\fk\ar[d]^{\mathrm{ev}}\\
\Gln{n}{\Bbbk}-\mathrm{mod}\ar[r]_{\mathcal{O}}&\gln{n}{\Bbbk}-\mathrm{mod}
}}

If $F,G\in\p_{\Bbbk}$ are homogeneous of degree less than $n$, then it follows from \cite{FS97} that the evaluation on $\Bbbk^{n}$ induces a natural isomorphism
\ali{\ext{\p_{\Bbbk}}{*}{F}{G}\to \ext{\Gln{n}{\Bbbk}}{*}{F\moc{\Bbbk^{n}}}{G\moc{\Bbbk^{n}}}.}
As the forgetful functors shown in the diagram \eqref{diagram3} are exact, we obtain the commutative square
\alilabel{comsquare}{\xymatrix{
\ext{\p_{\Bbbk}}{*}{F}{G}\ar[d]_{}\ar[r]^{}&\ext{\fk}{*}{F}{G}\ar[d]^{}\\
\ext{\gln{n}{\Bbbk}}{*}{F\moc{\Bbbk^{n}}}{G\moc{\Bbbk^{n}}}\ar[r]_{}&\ext{\Gln{n}{\Bbbk}}{*}{F\moc{\Bbbk^{n}}}{G\moc{\Bbbk^{n}}}
}}
According to the work of Cline, Parshall, Scott, and van der Kallen \cite{CPSvdK77}, Ext-groups of $\gln{n}{\Bbbk}-$modules can be computed via that of $\Gln{n}{\Bbbk}-$modules, with the help of the Frobenius twist. The square \eqref{comsquare} allows to establish similar relations between  Ext-groups in $\fk$ and that of strict polynomial functors. Let recall briefly about this fact. 

Given an integer $r\geq 1$, we denote by $I^{(r)}$ the strict polynomial functor defined as the intersection of the kernel of  $\p_{\Bbbk}-$morphisms 
$S^{p^{r}}\to S^{k}\otimes S^{p^{r}-k}$ for all $0<k<p^{r},$ 
induced by the comultiplication of the graded Hopf algebra $S^{*}$ (to be recalled in  Subsection \ref{twistpol}). The nontrivial strict polynomial functor $I^{(r)}$ is homogeneous of degree $p^{r}$. Therefore, $I^{(r)}$ and $I^{(s)}$ are not isomorphic if $r\neq s$. However, if $\Bbbk$ is perfect, then the forgetful functor  sends $I^{(r)}$ to the inclusion functor $I$ for all $r\geq 1$. The Frobenius twist of a strict polynomial functor $F$, denoted by $F^{(r)}$, is defined as $F\circ I^{(r)}$. The relation between Ext-groups in $\fk$ and in $\p_{\Bbbk}$ now goes as follows.
\begin{thm}[\cite{FFSS99}]
Let $\Bbbk$ be a perfect field, and let  $F,G\in\p_{\Bbbk}$ be two homogeneous strict polynomial functors of degree $d$, which is less than the cardinal of $\Bbbk$. If $r$ is big enough with respect to $i$, then the natural morphism 
\ali{\ext{\p_{\Bbbk}}{i}{F^{(r)}}{G^{(r)}}\to\ext{\fk}{i}{F}{G}}
is an isomorphism.
\end{thm}
In particular, we have the following corollary.
\begin{cor}\label{stabmclane}
If $r$ is big enough with respect to $i$, then the natural morphism
\ali{\ext{\p_{\f2}}{i}{I^{(r)}}{I^{(r)}}\to\ext{\mathscr{F}_{\f2}}{i}{I}{I}}
is an isomorphism.
\end{cor}
Therefore, we can compute the Mac Lane cohomology $\hml{\f2}{I}$ \textit{via} $\ext{\p_{\f2}}{i}{I^{(r)}}{I^{(r)}}$. In \cite{FS97}, Friedlander and Suslin compute these groups by constructing explicit injective resolution for $I^{(r)}$, using the method pioneered by Franjou, Lannes, and Schwartz \cite{FLS94}. We now show that the \phr method also allows to construct explicit resolutions of $I^{(r)}$ and leads to the same result.
\subsubsection{Injective resolutions of twisted strict polynomial functors}\label{twistpol} 
We begin with some basic homological algebra of strict polynomial functors. Throughout this subsection, we restrict our attention to the case $\Bbbk=\f2$, and we shorten $\p_{\f2,n}$ to $\p_{n}$.

\begin{defi}
An \textit{exponential strict polynomial functor} is a graded functor $F^{*}=\moc{F^{0},F^{1},F^{2},\ldots,F^{n},\ldots}$, where $F^{n}\in \p_{n}$, together with natural isomorphisms
\ali{F^{0}(V)=\f2, \quad F^{n}\moc{V\oplus W}\cong \bigoplus_{i=0}^{n}F^{i}\moc{V}\otimes F^{n-i}(W)}
for all integers $n>0$ and all $\f2-$vector spaces $V,W$.
\end{defi} 
Typical examples of exponential functors are given by the symmetric algebra $S^{*}$, the exterior algebra $\Lambda^{*}$, and the divided power algebra $\Gamma^{*}$. For all exponential strict polynomial functor $F^{*}$ and all $\f2-$vector space $V$, the addition map $+: V\oplus V\to V$ and the diagonal map $\Delta:V\to V\oplus V$ give rise to the natural maps
\ali{F^{n}(V)\otimes F^{m}(V)\hookrightarrow F^{n+m}(V\oplus V)\xrightarrow{F^{n+m}(+)} F^{n+m}(V),\\
F^{n+m}(V)\xrightarrow{F^{n+m}(\Delta)}F^{n+m}(V\oplus V)\twoheadrightarrow F^{n}(V)\otimes F^{m}(V),}
and then define natural product and coproduct operations 
\ali{\mathrm{mult}:F^{n}\otimes F^{m}\to F^{n+m},\quad \mathrm{comult}:F^{n+m}\to  F^{n}\otimes F^{m}.}
Remark that $S^{1}=\Lambda^{1}=\Gamma^{1}=I$. Let $F^{*},G^{*}\in\left\{S^{*},\Lambda^{*}\Gamma^{*} \right\}$, then we define $\kappa_{F^{*},G^{*}}^{d,e}$ as the composition map
\ali{F^{d}\otimes G^{e}\xrightarrow{\mathrm{comult}\otimes Id}F^{d-1}\otimes F^{1}\otimes G^{e}=F^{d-1}\otimes G^{1}\otimes G^{e}\xrightarrow{Id\otimes \mathrm{mult}}F^{d-1}\otimes G^{e+1},}
Following \cite{FLS94}, we have two exact sequences:
\begin{thm}\label{Koszulcomplexes}
The complexes $\moc{\left.\Gamma^{n}\otimes \Lambda^{d-n},\kappa_{\Gamma^{*},\Lambda^{*}}^{n,d-n}\right|0\leq n\leq d}$ and $\moc{\left.\Lambda^{n}\otimes S^{d-n},\kappa_{\Lambda^{*},S^{*}}^{n,d-n}\right|0\leq n\leq d}$ are exact. 
\end{thm}
Let $\lambda=\moc{\lambda_{1},\lambda_{2},\ldots,\lambda_{m}}\in \mathbb{N}^{m}$, and let  $F^{*}=\moc{F^{0},F^{1},F^{2},\ldots,F^{n},\ldots}$, where $F^{n}\in \p_{n}$, be an exponential strict polynomial functor. We denote by $F^{\lambda}$ the tensor product
\ali{F^{\lambda_{1}}\otimes F^{\lambda_{2}}\otimes\cdots \otimes F^{\lambda_{m}}.} Recall that the symmetric power functors $S^{d},d\geq 0,$ are injective in $\p_{\f2}$, and the tensor product of injective strict polynomial functors remains injective.  Moreover, the set $\left\{\left. S^{\lambda}\right| \lambda\in\mathbb{N}^{m},m\geq 1 \right\}$ forms a system of injective cogenerators of $\p_{\f2}$. 
According to \cite{FLS94}, the Frobenius twist of $S^{d}$ admits the following   explicit injective resolution.
\begin{thm}[{\cite[Section 2]{FLS94}}]
For all integers $d\geq 1,$ the complex 
$\moc{\left.S^{n}\otimes S^{d-n},\kappa_{S^{*},S^{*}}^{n,d-n}\right|0\leq d\leq n}$
is an injective resolution of $ S^{d(1)}$. Throughout this article, we denote this complex by $\mathscr{S}(d,1)$.
\end{thm}
Because the tensor product of injective strict polynomial functors is again injective, then $S^{\lambda(1)}$ admits \ali{\mathscr{S}(\lambda,1):=\mathscr{S}(\lambda_{1},1)\otimes \mathscr{S}(\lambda_{2},1)\otimes\cdots\otimes \mathscr{S}(\lambda_{m},1)} 
as an injective resolution. Since the Frobenius twist is exact, the \phr method is a suitable way to construct injective resolutions of $S^{\lambda(r)}$ for all integers $r\geq 1$. 

Let $\mathscr{S}$ denote the class of all resolutions that are direct sum of resolutions of the form $\mathscr{S}(\lambda,1)$. Hence, for all injective strict polynomial functors $J$, $J^{(1)}$ has a unique injective resolution belonging to $\mathscr{S}$. Given an injective resolution $J^{\bullet}$ of a strict polynomial functor $F$, in this paragraph, we denote by $\psr{F^{(1)}}{J^{\bullet(1)}}$ the \phr of $F^{(1)}$ with respect to the collection of resolutions of $J^{\bullet(1)}$ coming from $\mathscr{S}$. 

We define: 
\ali{\mathscr{S}(d,r):=\left\{\begin{array}{cl}
\mathscr{S}(d,1)&\textup{ if }r=1,\\
\psr{S^{d(r)}}{\mathscr{S}(d,r-1)^{(1)}}&\textup{ if }r>1.
\end{array} \right.}
The following lemma is straightforward.
\begin{lem}
The complex $\mathscr{S}(d,r)$ is an injective resolution of $S^{d(r)}$ for all integers $r\geq 1$. Moreover, the length of $\mathscr{S}(d,r)$ is $\moc{2^{r+1}-2}d$.
\end{lem}
We say that $\lambda=\moc{\lambda_{1},\lambda_{2},\ldots,\lambda_{m}}\in \mathbb{N}^{m}$ is divisible by $n$ if $\lambda_{k}$ is divisible by $n$ for all $1\leq k\leq m$. Remark that, if $\lambda$ is divisible by $2^{r}$, then the $s-$th term of the resolution $\mathscr{S}(\lambda,1)$ contains a factor of the form $S^{\alpha}$ such that $2^{r+1}|\alpha$ if and only if $2^{r+1}|s$. Otherwise, no factor of this form can appear in $\mathscr{S}(\lambda,1)$. A simple induction on $n$ using the \phr method shows that a similar property for $\mathscr{S}(d,n)$ also holds.
\begin{lem}
If $d$ is divisible by $2^{r}$, then the $s-$th term of the resolution $\mathscr{S}(d,n)$ contains a direct summand of the form $S^{\lambda}$ such that $2^{n+r}|\lambda$ if and only if $2^{n+r}|s$. Otherwise, no direct summand of this form can appear in $\mathscr{S}(d,n)$.  
\end{lem} 
In particular, when $d=1$, we obtain the following result.
\begin{cor}
The resolution $\mathscr{S}(1,n)$ is of length $2^{n+1}-2$. Denote by $\mathscr{S}(1,n)^{s}$ the $s-$th term of this resolution. Then, $\mathscr{S}(1,n)^{2k}$ contains a  unique direct summand of the form $S^{2^{n}}$ for all $0\leq k\leq 2^{n}-1$, and $\mathscr{S}(1,n)^{2k+1}$ contains no direct summand of this form.
\end{cor}
Recall that, for all $\lambda=\moc{\lambda_{1},\lambda_{2},\ldots,\lambda_{m}}\in \mathbb{N}^{m}$, we have
\ali{\ho{\p_{\f2}}{F^{(r)}}{S^{\lambda}}\cong \left\{\begin{array}{cl}
\ho{\p_{\f2}}{F}{S^{\moc{\frac{\lambda_{1}}{2^{r}},\frac{\lambda_{2}}{2^{r}},\ldots,\frac{\lambda_{m}}{2^{r}}}}}&\textup{ if }2^{r}|\lambda,\\
0&\textup{ otherwise.}
\end{array} \right.
}
Therefore, if $F$ is a homogeneous strict polynomial functor of degree $d$, then there is no nontrivial differentials in the complex $\ho{\p_{\f2}}{F^{(r)}}{\mathscr{S}(d,r)}$.
\begin{cor}
For all homogeneous strict polynomial functor $F$ of degree $d$, we have
\ali{\ho{\p_{\f2}}{F^{(r)}}{\mathscr{S}(d,r)}\cong\ext{\p_{\f2}}{*}{F^{(r)}}{S^{d(r)}}.}
In particular, we have
\ali{\ext{\p_{\f2}}{k}{I^{(r)}}{I^{(r)}}\cong \left\{\begin{array}{cl}
\f2&\text{ if }2|k, \text{ and }k\leq 2^{r+1}-2,\\
0&\text{ otherwise.}
\end{array} \right.
}
\end{cor}
Following Corollary \ref{stabmclane},  we obtain the main result of this subsection.
\begin{thm}[\cite{Bre78,FLS94,FS97}]\label{mclane}
We have:
\ali{\Hml{k}{\f2}{I}\cong \left\{\begin{array}{cl}
\f2&\text{ if }2|k,\\
0&\text{ otherwise.}
\end{array} \right.
}
\end{thm}

\subsection{Global dimension of homogeneous strict polynomial functors}
It is known that each homogeneous strict polynomial functor of finite degree is of finite projective (injective) dimension (see \cite{AB88,Don86,Tot97}). In this subsection, we will use the \phr method to show that the global dimension of the cateogry $\p_{d}$ is $2d-2$ for all $d\geq 0$.
\subsubsection{Resolutions of classical strict polynomial functors}
In this paragraph, we will construct explicit resolutions of the classical exponential strict polynomial functors $\Lambda^{*}$, $\Gamma^{*}$, and $S^{*}$. First, we fix the following notations. Let $F^{*}=\moc{F^{0},F^{1},F^{2},\ldots,F^{n},\ldots}$, where $F^{n}\in \p_{n}$, be an exponential strict polynomial functor. We denote:
\alilabel{recursion}{
A(d,s)&:=\left\{\moc{n_{1},n_{2},\ldots,n_{s}}\left| n_{i}\geq 1,\sum_{i=1}^{s}n_{i}=d  \right. \right\},\\
F(d,s)&:=\bigoplus_{A(d,s)}F^{\moc{n_{1},n_{2},\ldots,n_{s}}\label{recursion1}}.
}
Recall that the divided power functors $\Gamma^{d},d\geq 0,$ are projective in $\p_{\f2}$, and the tensor product of projective strict polynomial functors remains projective. Moreover, the set $\left\{\left. \Gamma^{\lambda}\right| \lambda\in\mathbb{N}^{m},m\geq 1 \right\}$ forms a system of projective generators of $\p_{\f2}$. We will use Theorem \ref{Koszulcomplexes} to construct explicit injective and projective resolutions of $\Lambda^{d}$. In fact, for $d=1$ and $d=2$, the complex $\moc{\left.\Gamma^{n}\otimes \Lambda^{d-n},\kappa_{\Gamma^{*},\Lambda^{*}}^{n,d-n}\right|0< n\leq d}$ is a projective resolution of $\Lambda^{d}$. In general, we construct a canonical projective resolution of $\Lambda^{d}$ as follows. We denote by $\Lambda_{d}$ the set of canonical projective resolution $C_{proj}\moc{\Lambda^{n}}$ of $\Lambda^{n}$ for all $n\leq d$, which is defined recursively as follows:
\ali{\Lambda_{2}&:=\left\{\left.\moc{\left.\Gamma^{n}\otimes \Lambda^{d-n},\kappa_{\Gamma^{*},\Lambda^{*}}^{n,d-n}\right|0< n\leq d}\right| 1\leq d\leq 2 \right\},\\
\Lambda_{d}&:=\Lambda_{d-1}\bigsqcup \left\{\psr{\Lambda^{d}}{\left\{\left.\Gamma^{d-s}\otimes C_{proj}\moc{\Lambda^{s}}\right| 0\leq s\leq d-1 \right\}}\right\}.
}
A simple induction on $d$ yields the following result.
\begin{lem}
For all integers $d\geq 1$, we have
\ali{\h\moc{t,C_{proj}\moc{\Lambda^{d}}}=\sum_{s=0}^{d-1}\Gamma(d,d-s)t^{s}.}
(See \eqref{recursion1} for the definition of $\Gamma(d,d-s)$.)
\end{lem}
Similarly, we construct a canonical injective resolution for $\Lambda^{d}$ as follows. We denote by $\mathrm{Inj}_{d}$ the set of canonical injective resolution $C_{inj}\moc{\Lambda^{n}}$ of $\Lambda^{n}$ for all $n\leq d$, which is defined recursively by:
\ali{\mathrm{Inj}_{2}&:=\left\{\left.\moc{\left.\Lambda^{n}\otimes S^{d-n},\kappa_{\Lambda^{*},S^{*}}^{n,d-n}\right|0\leq n\leq d}\right| 1\leq d\leq 2 \right\},\\
\mathrm{Inj}_{d}&:=\mathrm{Inj}_{d-1}\bigsqcup \left\{\psr{\Lambda^{d}}{\left\{\left.C_{inj}\moc{\Lambda^{s}}\otimes S^{d-s}\right| 0\leq s\leq d-1 \right\}}\right\}.
}
A simple induction on $d$ yields the following result.
\begin{lem}
For all $d\geq 1$, we have
\ali{\h\moc{t,C_{inj}\moc{\Lambda^{d}}}=\sum_{s=0}^{d-1}S(d,d-s)t^{s}.}
(See \eqref{recursion} for the definition of $S(d,d-s)$.)
\end{lem}
Now, using the complex $\moc{\left.\Gamma^{n}\otimes \Lambda^{d-n},\kappa_{\Gamma^{*},\Lambda^{*}}^{n,d-n}\right|0< n\leq d}$, the set $\mathrm{Inj}_{d}$, and the \phr method, we can define a set $\Gamma_{d}$ of canonical injective resolution $C_{inj}\moc{\Gamma^{n}}$ of $\Gamma^{n}$ for all $n\leq d$ as follows:
\ali{
\Gamma_{1}&:=\left\{\left.\moc{\left.\Gamma^{n}\otimes \Lambda^{d-n},\kappa_{\Gamma^{*},\Lambda^{*}}^{n,d-n}\right|0< n\leq d}\right| d=1 \right\},\\
\Gamma_{d}&:=\Gamma_{d-1}\bigsqcup \left\{\psr{\Gamma^{d}}{\left\{\left.C_{inj}\moc{\Gamma^{s}}\otimes C_{inj}\moc{\Lambda^{d-s}}\right| 0\leq s\leq d-1 \right\}}\right\}.}
Induction on $d$ gives rise to the following lemma.
\begin{lem}\label{injdim}
For all integers $d\geq 1$, we have
\ali{\h\moc{t,C_{inj}\moc{\Gamma^{d}}}=\sum_{s=0}^{d-1}S(d,d-s)t^{2s}(1+t)^{d-s-1}.}
\end{lem}
Similarly, using the complex $\moc{\left.\Lambda^{n}\otimes S^{d-n},\kappa_{\Lambda^{*},S^{*}}^{n,d-n}\right|0\leq n\leq d}$, the set $\Gamma_{d}$, and the \phr method, we can define a set $S_{d}$ of canonical projective resolution $C_{proj}\moc{S^{n}}$ of $S^{n}$ for all $n\leq d$ as follows:
\ali{
S_{1}&:=\left\{\left.\moc{\left.\Lambda^{n}\otimes S^{d-n},\kappa_{\Lambda^{*},S^{*}}^{n,d-n}\right|0< n\leq d}\right| d=1 \right\},\\
S_{d}&:=S_{d-1}\bigsqcup \left\{\psr{S^{d}}{\left\{\left.C_{proj}\moc{\Lambda^{d-s}}\otimes C_{proj}\moc{S^{s}}\right| 0\leq s\leq d-1 \right\}}\right\}.}
Induction on $d$ gives rise to the following lemma.

\begin{lem}\label{projdim}
For all $d\geq 1$, we have:
\ali{\h\moc{t,C_{proj}\moc{S^{d}}}=\sum_{s=0}^{d-1}\Gamma(d,d-s)t^{2s}(1+t)^{d-s-1}.}
\end{lem}

\begin{rem}\label{globdim}
Lemma \ref{injdim} and \ref{projdim} signify that 
\ali{\mathrm{dim}_{inj}\moc{\Gamma^{d}}\leq 2d-2,\quad \mathrm{dim}_{proj}\moc{S^{d}}\leq 2d-2.}
\end{rem}
\subsection{Global dimension of strict polynomial functors} 
The following lemma is a direct consequence of Corollary \ref{projdim}.
\begin{lem}\label{5.16}
For all integers $n\geq 1$, we have
\ali{\ext{\p_{2^{n}}}{2^{n+1}-2}{S^{2^{n}}}{\Gamma^{2^{n}}}\cong\f2.}
Hence, $S^{2^{n}}$ is of projective dimension $2^{n+1}-2$, and $\Gamma^{2^{n}}$ is of injective dimension $2^{n+1}-2$.
\end{lem}
\begin{proof}
Denote $2^{n}$ by $d$. In view of Corollary \ref{projdim}, $\ext{\p_{d}}{2d-2}{S^{d}}{\Gamma^{d}}$ is isomorphic to the cokernel of the map
\ali{f:\ho{\p_{d}}{\Gamma\moc{d,2}}{\Gamma^{d}}\to \ho{\p_{d}}{\Gamma^{d}}{\Gamma^{d}}.}
If this cokernel was trivial, then $\Gamma^{d}$ would be a direct summand of $\Gamma\moc{d,2}$, whence a contradiction. Therefore, $\cke{f}$ is nontrivial. But, on the other hand, $\ho{\p_{d}}{\Gamma^{d}}{\Gamma^{d}}\cong \f2$, then the lemma follows from Remark \ref{globdim}. 
\end{proof}
\begin{cor}\label{particulardim}
Let $d\geq 1$ be an integer and let $\sum\limits_{i=1}^{k}2^{n_{k}}$ be its $2-$adic expression. Let $I\in\p_{d}$ be an injective strict polynomial functors, and let $P\in\p_{d}$ be a projective one. Then, we have
\ali{\mathrm{dim}_{inj}P\leq 2d-2k,\quad \mathrm{dim}_{proj}I\leq 2d-2k.} 
In particular, let $\lambda=\moc{2^{n_{1}},2^{n_{2}},\ldots,2^{n_{k}}}$, then we have:
\ali{\mathrm{dim}_{inj}\Gamma^{\lambda}= 2d-2k,\quad \mathrm{dim}_{proj}S^{\lambda}= 2d-2k.}
\end{cor}
\begin{proof}
The corollary follows from Lemma \ref{5.16} and from the fact that the morphism $\Gamma^{\lambda}\to\Gamma^{d}$ induced by the multiplication of $\Gamma^{*}$ is a splitting projection.
\end{proof}
Corollary \ref{particulardim} can be extended to all strict polynomial functors.
\begin{lem}
Let $d\geq 1$ be an integer and let $\sum\limits_{i=1}^{k}2^{n_{k}}$ be its $2-$adic expression. Let $F\in\p_{d}$ be an arbitrary strict polynomial, then 
\ali{\mathrm{dim}_{inj}F\leq 2d-2k,\quad \mathrm{dim}_{proj}F\leq 2d-2k.} 
\end{lem}
\begin{proof}
It is sufficient to show that
\ali{\ext{\p_{d}}{s}{G}{F}&=0,\textup{ for all }s>2d-2 \textup{ and all }G\in\p_{d},\\
\ext{\p_{d}}{s}{F}{G}&=0,\textup{ for all }s>2d-2 \textup{ and all }G\in\p_{d}.}
But these computations can be carried out by induction on injective and projective dimensions of $F$, using Corollary \ref{particulardim}.
\end{proof}

\begin{thm}[\cite{Tot97}]\label{globaldimension}
Let $d\geq 1$ be an integer and let $\sum\limits_{i=1}^{k}2^{n_{k}}$ be its $2-$adic expression. Then, the global dimension of $\p_{d}$ is $2d-2k$.
\end{thm}

\section{The Lambda algebra}\label{Lambdaalgebra}
This section aims at giving a new construction of the Lambda algebra, using the \phr method.

\subsection{Iterated suspensions of Brown-Gitler modules}\label{iterBG}
For all integers $n\geq 0$, the Mahowald short exact sequence
$$0\to \Sigma J(n)\xrightarrow{s_{n}}J(n+1)\xrightarrow{\bullet Sq^{\frac{n+1}{2}}}J\left(\frac{n+1}{2}\right)\to 0$$
yields an injective resolution for $\Sigma J(n)$. Since the functor $\Sigma$ is exact, applying it to an injective resolution of  $\Sigma^{m} J(n)$ induces an acyclic sequence such that its unique nontrivial cohomology is $\Sigma^{m+1} J(n)$. Therefore, the \phr method allows to construct an injective resolution of $\Sigma^{m+1} J(n)$ from that of $\Sigma^{m} J(n)$.
\begin{defi}
Let $B=\bigoplus_{A}J(n_{\alpha})$ be a BG module, we define: \begin{align*}
\aug{}{B}&:=\bigoplus_{A}J(1+n_{\alpha}).
\end{align*} 
\end{defi}

A reformulation of Lemma \ref{key9} gives rise to the following proposition.
\begin{prop}[\Phr for the suspension of BG complexes]\label{BGS} Let
	$\left(B^{k},\tau^{k}:B^{k}\to B^{k+1},\kappa\geq 0\right)$ be a BG complex. Denote  by $p^{k}$ the canonical projection $\aug{}{B^{k}}\to\tphi{\aug{}{B^{k}}}$ induced by the natural transformation $Id\to\tilde{\Phi}$. Then, for all integers $k\geq 0,$ there exist morphisms
	\begin{align*}
	\theta^{k}: \aug{}{B^{k}}&\to  \aug{}{B^{k+1}},\\
	\delta^{k}: \tphi{\aug{}{B^{k}}}&\to \aug{}{B^{k+2}},\\
	\omega^{k}:\tphi{\aug{}{B^{k}}}&\to \tphi{\aug{}{B^{k+1}}} 
	\end{align*}
	such that the morphisms
	$\left(\begin{smallmatrix}
	\theta^k &\delta^{k-1}\\
	p^{k} & \omega^{k-1}
	\end{smallmatrix}\right)$ for all $k\geq 1,$ denoted by $\partial^{k}$, and $\left(\begin{smallmatrix}
	\theta^0\\
	p^{0}
	\end{smallmatrix}\right)$, denoted by $\partial^{0},$
	make the sequence 
	\begin{equation}
	0\to \aug{}{B^{0}}\xrightarrow{\partial^{0}} \aug{}{B^{1}}\bigoplus
	\tphi{\aug{}{B^{0}}}\xrightarrow{\partial^{1}}\cdots\xrightarrow{\partial^{k}}\aug{}{B^{k+1}}\bigoplus
	\tphi{\aug{}{B^{k}}}\xrightarrow{\partial^{k+1}}\cdots
	\end{equation}
	a complex that we denote by $\pb{}{B^{k},\tau^{k}}$. It is acyclic provided that $\left(B^{k},\tau^{k}\right)$ is acyclic either, and in this case, the only nontrivial cohomology of $\pb{}{B^{k},\tau^{k}}$ is the suspension of that of $\left(B^{k},\tau^{k}\right)$. 
\end{prop}
\begin{cor}\label{48}
	Denote by $J_{n}^{\bullet}$ the acyclic complex where $J_{n}^{0}=J(n)$ and $J_{n}^{k}$ is trivial for all integers $k\geq 1$. Then,  $\pb{m}{J_{n}^{\bullet}}$ is an injective resolution of $\Sigma^{m}J(n)$. Hereafter, we denote this resolution by $\pb{}{m,n}$.
\end{cor}

\begin{defi}\label{gmn}
For all integers $m,n\geq 0$, each term of $\pb{}{m,n}$ is a direct sum of Brown-Gitler modules. Denote by $G(m,n)$ the graph associated with $\pb{}{m,n}$ with respect to this decomposition, and $V_{r}\moc{G(m,n)},E_{r}\moc{G(m,n)}$ are shortened to $V_{r}(m,n),E_{r}(m,n)$ respectively. (See Section \ref{graphrep} for the definition of graph representation.)
\end{defi}

\begin{depo}\label{BGS1}
For all integers $m,n\geq 0$, we have:
\begin{enumerate}
\item\label{item1} There exists a canonical inclusion of graphs $i_{m,n}:G(m,n)\to G(m+1,n)$ corresponding to the inclusion of complexes $\Sigma \pb{}{m,n}\to \pb{}{m+1,n}$.
\item\label{item2} The set $V_{0}(m,n)$ contains a unique vertex, denoted by $p_{m,n}$, corresponding to $J(m+n)$.
\item\label{item3} For all integers $r\geq 0$ and all $v\in V_{r}(m,n),w\in V_{r+1}(m,n)$, then either $[v,w]=\bullet Sq^{k}$ for some integer $k\geq 0$, or $[v,w]=0$. Hereafter, for the sake of simplicity, we write  $[v,w]=Sq^{k}$ instead of $[v,w]=\bullet Sq^{k}$.
\item\label{item4} For all vertices $v\in V_{r}(m,n)$, there exists a unique collection of vertices 
\ali{\left\{ v_{i}\left|\begin{array}{l}
v_{i}\in V_{i}(m,n), 0\leq i\leq r,\\
v_{0}=p_{m+n}, v_{r}=v,\\
\left[v_{i},v_{i+1}\right]=Sq^{k_{i}},\\
\left[v_{i},w\right]\neq Sq^{k_{i}},\forall v_{i+1}\neq w,\\ 
2k_{0}> n,\\
2k_{i+1}>k_{i}, r-2\geq i\geq 0,\\
m+n-\sum_{i=0}^{r-1}k_{i}\geq k_{r-1},
\end{array}\right.\right\}}
forming a path from $p_{m+n}$ to $v$, and such a path is called an \textit{admissible path}. In this case, we say that $v$ is of bidegree $\moc{r,\sum_{i=0}^{r-1}k_{i}}$. We call this the \textit{Lambda bidegree} of $v$ and denote it by $\norm{v}=\moc{\norm{v}_{1},\norm{v}_{2}}$, where $\norm{v}_{1}=r$ and $\norm{v}_{2}=\sum_{i=0}^{r-1}k_{i}$.
\item\label{item5} The set $E(m,n)$ of edges of $G(m,n)$ can be defined as the colimit of the increasing sequence of sets 
\ali{E^{0}(m,n)\subset E^{1}(m,n)\subset\cdots\subset E^{r}(m,n)\subset\ldots}
as follows. Denote by $E^{0}(m,n)$ the set of all edges induced by the admissible paths. Suppose that we have defined $E^{r}(m,n)$. We now define $E^{r+1}(m,n)$.  For all vertices $v\in V_{i}(m,n)$ and $w\in V_{i+2}(m,n)$, consider the sum of all products $\mocv{v,u}\mocv{u,w}$ over all $u\in V_{i+1}(m,n)$, where $\mocv{v,u},\mocv{u,w}$ belong to $E^{r}(m,n)$. Recall that this is a sum of products of Steenrod operations of the form $Sq^{k}$. Then, we can write this sum in terms of admissible monomials $\sum_{\moc{i_{1},i_{2}}\in A}Sq^{i_{1}}Sq^{i_{2}}$. Let $z\in V_{i+1}(m,n)$ such that $\mocv{v,z}=Sq^{i_{1}}$ belongs to $E^{0}(m,n)$, then we connect $z$ to $w$ by the edge $\mocv{z,w}=Sq^{i_{2}}$. We now define $E^{r+1}(m,n)$ as the union of $E^{r}(m,n)$ with these new edges. 
\end{enumerate}
\end{depo}
\begin{proof}
The proof is carried out by induction on $m$. The case $m=1$ is straightforward as $G(1,n)$ is the graph corresponding to the Mahowald short exact sequence \eqref{mahowalds} of $\Sigma J(n)$. Suppose that the proposition is proved for all $m<k$, where $k\geq 2$ is a certain integer. We now verify the case $m=k$. By definition, we have
\ali{ \pb{}{k,n}^{s}\cong \aug{}{ \pb{}{k-1,n}}^{s}\oplus \tilde{\Phi}\moc{\aug{}{ \pb{}{k-1,n}}}^{s-1}}
for all integers $s\geq 0$, with the convention that $\tilde{\Phi}\moc{\aug{}{ \pb{}{k-1,n}}}^{-1}=0$. Here, the upper index $s$ signifies the $s-$th term of the complex. Therefore, $\pb{}{k,n}^{0}=\aug{}{ \pb{}{k-1,n}}^{0}$, and hence, following the induction hypothesis, the set $V_{0}(k,n)$ contains a unique vertex, denoted by $p_{k,n}$, corresponding to $J(k+n)$. \textbf{This verifies the point \eqref{item2}}. Moreover, the morphisms
\ali{\Sigma\pb{}{k-1,n}^{s}\to \aug{}{ \pb{}{k-1,n}}^{s}\hookrightarrow \pb{}{k,n}^{s},}
for all $s\geq 0$, induce a canonical inclusion $\Sigma \pb{}{k-1,n}\to \pb{}{k,n}$. As the graph of  $\Sigma \pb{}{k-1,n}$ and that of $\pb{}{k-1,n}$ are identical, we get an inclusion $G(k-1,n)\hookrightarrow G(k,n)$. \textbf{This concludes the proof of the point \eqref{item1}}.

Remark that, in view of Proposition \ref{BGS}, if $\tau^{s}$ is represented by a matrix with coefficients of the form $\bullet Sq^{s}$, then so are $\theta^{s},\delta^{s}$, and $\omega^{s}$. Therefore, each edge of $E_{r}(k,n)$ is also labeled by a Steenrod operation of the form $Sq^{s}$ for some integer $s$. \textbf{Then, the point \eqref{item3} holds}.

We also observe that the vertices of $G(k,n)$ that do not belong to the image of $G(k-1,n)$ correspond to $\tilde{\Phi}\moc{\aug{}{ \pb{}{k-1,n}}}$. \textbf{Then, the induction hypothesis and the Mahowald short exact sequences take care of the points \eqref{item4} and \eqref{item5}}.
\end{proof}

\begin{rem}\label{rmk66}
\begin{enumerate}
\item The bigrading structure of $G(m,n)$ is different from that of the general associated graph of a BG complex (see Definition \ref{graphdeg}). We make such a modification so that the canonical inclusion $i_{m,n}$ from $G(m,n)$ to  $G(m+1,n)$ is of bidegree $(0,0)$.
\item\label{rmk662} Let $k\geq 1$ be an integer. Then, for all integers $m,n$ such that $m+n\geq 2k$, the inclusion $i_{m,n}$ of $G(m,n)$ into $G(m+1,n)$ is an isomorphism on the area of bidegrees $(r,s)$ for all $s\leq k$. It is obvious that this is an isomorphism on the set of vertices. What is less evident is that it is also bijective on the set of edges. But, this is a consequence of the instability condition: the morphism $\bullet Sq^{q}:J(t+q)\to J(t)$ is nontrivial if and only if $q\leq t$. In fact, Proposition \ref{BGS} shows that all the edges of $G(m+1,n)$ that are not in the image of $i_{m,n}$ must have a vertex belonging to $\tilde{\Phi}\moc{\aug{}{ \pb{}{m,n}}}$, which is not in the area of consideration.
\end{enumerate}
\end{rem}
Therefore, using Lemma \ref{lambdaext} with an appropriate changing of bidegrees, we obtain the following stabilization result:
\begin{prop}\label{stabilization}
Let $k\geq 1$ be an integer. Then, for all for all integers $m,n,s$ such that $m+n\geq 2k$, and that $s\geq m+n- k$, the morphism
\ali{\ext{\u}{*}{\Sigma^{m}\f2}{\Sigma^{m+s}J(n)}\to \ext{\u}{*}{\Sigma^{m+1}\f2}{\Sigma^{m+s+1}J(n)}}
induced by the suspension $\Sigma$ is an isomorphism.
\end{prop}

\begin{defi}
For all integers $n\geq 1$, we denote by $G(n)$ the limit of the system
\ali{\cdots\xrightarrow{i_{m-1,n}}G(m,n)\xrightarrow{i_{m,n}}G(m+1,n)\xrightarrow{i_{m+1,n}}G(m+2,n)\xrightarrow{i_{m+2,n}}\cdots .}
\end{defi}
The following description of $G(n)$ follows from Proposition \ref{BGS1}.
\begin{lem}\label{gn}
For all integers $n\geq 0$, the set of vertices of $G(n)$ is the disjoint union of $V_{r}(n)$, and the set of edges is the disjoint union of $E_{r}(n)$, where $r\geq 0$, such that:
\begin{enumerate}
\item The set $V_{0}(n)$ contains a unique vertex, denoted by $p_{n}$.
\item For all integers $r\geq 1$, each edge of $E_{r}(n)$, from $v\in V_{r}(n)$ to $w\in V_{r+1}(n)$, corresponds to a morphism of the form $\bullet Sq^{k}$ for some integer $k\geq 0$, and we label this edge by $Sq^{k}$. In this case, we write $[v,w]=Sq^{k}$. If there is no edge from $v$ to $w$, then we write $[v,w]=0$. 
\item For all vertices $v\in V_{r}(n)$, there exists a unique collection of vertices 
\ali{\left\{ v_{i}\left|\begin{array}{l}
v_{i}\in V_{i}(n), 0\leq i\leq r,\\
v_{0}=p_{m+n}, v_{r}=v,\\
\left[v_{i},v_{i+1}\right]=Sq^{k_{i}},\\
\left[v_{i},w\right]\neq Sq^{k_{i}},\forall w\neq v_{i+1},\\
2k_{0}> n,\\
2k_{i+1}>k_{i}, 0\leq i\leq r-2,\\
m+n-\sum_{i=0}^{r-1}k_{i}\geq k_{r-1},
\end{array}\right.\right\}}
forming a path from $p_{m+n}$ to $v$, and such a path is called an \textit{admissible path}. We also say that $v$ is of Lambda bidegree $\norm{v}=\moc{\norm{v}_{1},\norm{v}_{2}}=\moc{r,\sum_{i=0}^{r-1}k_{i}}$.
\end{enumerate}
In particular, $V_{1}(1)$ consists of vertices $\left\{\left.\lambda_{i} \right| 0\leq i, \left[p_{1},\lambda_{i} \right]=Sq^{i+1}\right\}$.
\end{lem}
\begin{rem}
In Lemma \ref{gn}, we do not describe the set of edges because it can be obtained by taking the colimit of the sets $E(m,n)$. And thanks to the point \eqref{rmk662} of Remark \ref{rmk66}, in practice, the computations concerning $G(n)$ will always be carried out using the model $G(m,n)$ for some suitable integer $m\geq 0$.
\end{rem}
\subsection{The Lambda algebra}
In this paragraph, we will show that there exists an appropriate product $\times$ on $G(1)$ such that $\moc{G(1),\times}$ is the Lambda algebra.
\begin{defi}
We define $\mathbb{T}$ to be the free bigraded algebra over $\f2$ generated by the symbols $\lambda_{i}$ of bidegree $\moc{1,i+1}$ for all integers $i\geq 0$, and $\Lambda$ to be the bigraded $\f2-$vector space generated by the set of vertices of $G(1)$, where we take into account the Lambda bidegree as defined in Lemma \ref{gn}.
\end{defi}

\begin{lem}
The $\f2-$vector space $\Lambda$ is endowed with a structure of right $\mathbb{T}-$module, defined as follows:
\alilabel{module}{x\lambda_{i}=\sum_{\left[x,y\right]=Sq^{i+1}}^{}y}
for all vertices $x$ of $G(1)$. Moreover, the $\mathbb{T}-$linear morphism
\ali{g:\mathbb{T}&\longrightarrow \Lambda\\
\alpha&\longmapsto p_{1}\alpha}
is an epimorphism, and the set
\ali{\left\{\left. p_{1}\lambda_{k_{1}}\lambda_{k_{2}}\ldots \lambda_{k_{n}}\right| 2k_{i}\leq k_{i+1},1\leq i\leq n-1 \right\}}
forms an $\f2-$basis of $\Lambda$.
\end{lem}
\begin{proof}
It is straightforward that the action law \eqref{module} defines the structure of a right $\mathbb{T}-$module for $\Lambda$. We now show that $g$ is surjective. Indeed, for all vertices $v\in V_{r}(1)$, there exists a unique admissible path \ali{\left\{ v_{i}\left|\begin{array}{l}
v_{i}\in V_{i}(1), 0\leq i\leq r,\\
v_{0}=p_{1}, v_{r}=v,\\
\left[v_{i},v_{i+1}\right]=Sq^{k_{i}},\\
\left[v_{i},w\right]\neq Sq^{k_{i}},\forall w\neq v_{i+1},\\
2k_{0}> n,\\
2k_{i+1}>k_{i}, r-2\geq i\geq 0
\end{array}\right.\right\}} connecting $p_{1}$ and $v$. Hence, we have
\ali{v=p_{1}\lambda_{k_{0}-1}\lambda_{k_{1}-1}\lambda_{k_{2}-1}\ldots \lambda_{k_{r-1}-1}.}
Then, this concludes that $g$ is surjective.
\end{proof}
\begin{lem}\label{lem613}
The kernel of $g:\mathbb{T}\to\Lambda$ is the two-sided ideal generated by:
\alilabel{twosidedideal}{\lambda_{a}\lambda_{b}=\sum_{2a+2b>3i+1\geq 6b+4}^{}\binom{a-2i-1}{i-2b-1}\lambda_{i}\lambda_{a+b-i}}
for all integers $a,b$ such that $a\geq 2b+1$.
\end{lem}
\begin{proof}
Recall that the set
\ali{\left\{\left. p_{1}\lambda_{k_{1}}\lambda_{k_{2}}\ldots \lambda_{k_{n}}\right| 2k_{i}\leq k_{i+1},1\leq i\leq n-1 \right\}}
forms an $\f2-$basis of $\Lambda$. We now show how to express $p_{1}\lambda_{a}\lambda_{b}$ as a linear combination of this basis for all $a\geq 2b+1$. Let $A$ be the set
\ali{\left\{ v\in V_{2}(1)\left| \left[p_{1}\lambda_{a},v\right]=Sq^{b+1}\right. \right\},}
then 
\ali{p_{1}\lambda_{a}\lambda_{b}=\sum_{v\in A}^{}v.}
By definition of $G(1)$ as the colimit of $G(m,1),m\geq 1,$ then $p_{1}\lambda_{a}\in G(m,1)$ if and only if $m\geq 2a+1$. Moreover, for all $w\in V_{2}(1)$ such that $[p_{1}\lambda_{a},w]\neq 0,$ then $w\in V_{2}(m,1)$ for all $m\geq 2a+1$, and $[p_{1}\lambda_{a},w]\neq 0$ in $G(m,1)$. In particular, when $m=2a+1$, then $p_{1}\lambda_{a}$ corresponds to a direct summand of the form $J(a+1)$ of $\tilde{\Phi}\moc{\aug{}{ \pb{}{k-1,n}}}^{0}$. Therefore, an element $w\in V_{2}(2a+1,1)$, with $\left[p_{1}\lambda_{a},w\right]=Sq^{b+1},$ where $a\geq 2b+1,$ exists if and only if $w$ is defined by the admissible form $p_{1}\lambda_{i}\lambda_{j}$ such that $Sq^{a+1}Sq^{b+1}$ appears in the expression of $Sq^{i+1}Sq^{j+1}$ as linear combination of admissible Steenrod monomials. Hence, it follows from the Adem relations that we have 
\ali{p_{1}\lambda_{a}\lambda_{b}=\sum_{2a+2b>3i+1\geq 6b+4}^{}\binom{a-2i-1}{i-2b-1}p_{1}\lambda_{i}\lambda_{a+b-i}}
for all $a\geq 2b+1$. Now, let $x=p_{1}\lambda_{a_{1}}\lambda_{a_{2}}\ldots \lambda_{a_{k}}$, and denote by $m<k$ the first index such that $a_{m}\geq 2 a_{m+1}+1$. Then, by the same method, we have
\ali{p_{1}\lambda_{a_{1}}\lambda_{a_{2}}\ldots \lambda_{a_{k}}=p_{1}\lambda_{a_{1}}\lambda_{a_{2}}\ldots \lambda_{a_{m-1}}\moc{\sum_{3i= 6a_{m+1}+3}^{2a_{m}+2a_{m+1}-2}\binom{a_{m}-2i-1}{i-2a_{m+1}-1}\lambda_{i}\lambda_{a+b-i}}\lambda_{a_{m+2}}\lambda_{a_{m+3}}\ldots \lambda_{a_{k}}.}
As a result, the kernel of $g$ is the two-sided ideal generated by the relations \eqref{twosidedideal}.
\end{proof}
It follows from Lemma \ref{lambdaext} and Proposition \ref{stabilization} that $\Lambda$ is endowed with a differential graded module structure with explicit homology.
\begin{depo}\label{lambdaalgebra}
Let $d:\Lambda\to\Lambda$ be the $\f2-$linear map defined by:
\ali{d(x)=\sum_{\left[x,v \right]=Sq^{0}}^{}v}
for all vertices $x\in V_{r}(1)$. Then, $\moc{\Lambda,d}$ is a differential graded $\f2-$vector space. Moreover, we have
\ali{\h^{r,s}\moc{\Lambda,d}\cong \varinjlim\limits_{n}\ext{\u}{r}{\Sigma^{n}\f2}{\Sigma^{n+s}\f2}.}
\end{depo}
We will actually prove that $\moc{\Lambda,d}$ is a differential graded algebra.
\begin{defi}
Denote by $X_{}^{\bullet}$ the exact sequence where $X_{}^{0}=X_{}^{1}=J(1)$, $\partial^{0}=\bullet Sq^{0}$ and $X_{}^{k}$ is trivial for all $k\geq 2$. Then $\pb{m}{X_{}^{\bullet}}$ is an exact sequence that we denote by $X_{m}$. We denote by $G\moc{X_{m}}$ the associated graph with respect to the decomposition of the terms of $X_{m}$ as direct sum of Brown-Gitler modules.
\end{defi}
\begin{rem}
\begin{enumerate}
\item For all integers $m\geq 1$, we have
\ali{V_{r}\moc{G\moc{X_{m}}}&=V_{r}(m,1)\bigsqcup V_{r-1}(m,1),\\
E_{r}\moc{G\moc{X_{m}}}&=E_{r,r+1}(m,1)\bigsqcup E_{r-1}(m,1)\bigsqcup F_{r},}
where
\ali{F_{r}:=\left\{\left.[v,v]=Sq^{0}\right| v\in V_{r}(m,1)\subset V_{r}\moc{G\moc{X_{m}}} \right\}.}
\item For all integers $m\geq 1$, we have
\ali{\h\moc{t,G\moc{X_{m}}}=(1+t)\h\moc{t,G(m,1)}.}
\end{enumerate}
\end{rem}
Similar to the construction of $\Lambda$, we have the following result.
\begin{prop}
There exist canonical inclusions of graphs $G(X_{m})\to G(X_{m+1})$ corresponding the the inclusion of complexes $\Sigma X_{m}\to X_{m+1}$. Denote by $X$ the colimit of 
\ali{\cdots\to G(X_{m-1})\to G(X_{m})\to G(X_{m+1})\to\cdots,}
and by $\mathbb{L}$ the $\f2-$vector space generated by the vertices of $X$. Then, $\mathbb{L}$ is an associative bigraded algebra generated by the symbols $\lambda_{i}$ of bidegree $(1,i+1)$ for all integers $i\geq -1$, subject to the relations
\ali{\lambda_{a}\lambda_{b}=\sum_{2a+2b>3i+1\geq 6b+4}^{}\binom{a-2i-1}{i-2b-1}\lambda_{i}\lambda_{a+b-i},}
for all integers $a,b$ such that $a\geq 2b+1$.
\end{prop}
We now show that $\moc{\Lambda,d}$ is a differential algebra.
\begin{lem}
The algebra $\Lambda$ is a subalgebra of $\mathbb{L}$. Moreover, we have
\ali{d\moc{\alpha}=\alpha \lambda_{-1}-\lambda_{-1}\alpha.}
for all $\alpha\in\Lambda\subset \mathbb{L}$. Hence, $\moc{\Lambda,d}$ is a differential algebra.
\end{lem}
\begin{proof}
For all $\alpha\in\Lambda$ of bidegree $(r,s)$, we have:
\ali{\alpha\lambda_{-1}&=\sum_{\left[\alpha,v\right]=Sq^{0}}^{}v\\
&=\moc{\sum_{\substack{\left[\alpha,v\right]=Sq^{0}\\v\in V_{r+1}(1)}}v}+\moc{\sum_{\substack{\left[\alpha,v\right]=Sq^{0}\\v\in V_{r}(1)\subset V_{r+1}(X)}}v}\\
&=d(\alpha)+\lambda_{-1}\alpha.}
The lemma follows.
\end{proof}

\begin{defi}\label{lambdamn}
For all integers $m,n\geq 0$, let $\Lambda(m,n)$ denote the bigraded $\f2-$vector space generated by the vertices $V(m,n)$ of the graph $G(m,n)$, where the bidegree of the generator associated with a vertex $v$ of bidegree $(r,s)$ in $G(m,n)$ is $(r,m+n-s)$.
\end{defi}
\begin{rem}
\begin{enumerate}
\item It is obvious that the law
\alilabel{module}{x\lambda_{i}=\sum_{\left[x,y\right]=Sq^{i+1}}^{}y}
for all vertices $x$ of $G(m,n)$ defines an action of $\mathbb{T}$ on $\Lambda(m,n)$. Moreover, this action is unstable in the following sense: $x\lambda_{i}=0$ for all generators $x$ of bidegree $(r,s)$ such that $2i+2>s$. Using a similar argument as that in the proof of Lemma \ref{lem613}, we can show that this action passes to the quotient $\Lambda$ and yields a structure of right $\Lambda-$module for $\Lambda(m,n)$. 
\item We can also define a $\f2-$linear map
$\partial^{m,n}:\Lambda(m,n)\to\Lambda(m,n)$ by:
\ali{\partial^{m,n}(x)=\sum_{\left[x,v \right]=Sq^{0}}^{}v}
for all vertices $x\in V(m,n)$. In fact, we have
\ali{\partial^{m,n}\moc{p_{m,n}\theta}=p_{m,n}d\moc{\theta},}
where $d$ is the differential of $\Lambda$.
\end{enumerate}
\end{rem}

\section{Minimal resolutions and finite exact sequences of BG modules}\label{minimalres}
The Bockstein sequence, as it is defined in Definition  \ref{bocksteins}, plays an important role in the rest of the present paper. This section is devoted to study the interaction between the \phr method and the iterated suspensions of the Bockstein sequence. In fact, we will prove that this sequence is \textit{saturated} with respect to the \phr method, in the sense of Proposition \ref{saturation}.

Let $\mathcal{C}$ be an abelian category with enough injective objects, and let $M$ be an object of $\mathcal{C}$. Recall that an injective resolution $\moc{I^{i},\partial^{i},i\geq 0}$ of $M$ is minimal if $I^{0}$ is the injective hull of $M$, $I^{1}$ is the injective hull of $I^{0}/M$, and $I^{j}$ is that of $I^{j-1}/\im{\partial^{j-2}}$ for all $j\geq 2$. The following lemma is classical and is left to the readers.
\begin{lem}\label{injhull}
	Let $M$ and $N$ be two objects of an abelian category $\mathcal{C}$. Moreover,  let $\moc{I^{i},\partial^{i},i\geq 0}$ be a minimal injective resolution of $M$, and let $\moc{J^{i},\delta^{i},i\geq 0}$ be an injective resolution of $M\oplus N$. Then a split projection $M\oplus N\to M$ gives rise to a split projection of complexes 
	$\left\{\alpha^{i}:J^{i}\to I^{i},i\geq 0\right\}.$
\end{lem}
As a consequence, minimal injective resolutions are unique up to isomorphisms. 

Let denote by $\mathrm{Ch}(\mathcal{C})$ the category of cochain complexes $\moc{C^{i},\partial^{i},i\geq 0}$ of objects in $\mathcal{C}$. 
\begin{defi}[Extension of complexes]
If 
$$0\to \moc{C^{i},\partial^{i},i\geq 0}\xrightarrow{\alpha^{i}}\moc{D^{i},\delta^{i},i\geq 0}\xrightarrow{\beta^{i}}\moc{E^{i},\eta^{i},i\geq 0}\to 0$$ 
is a short exact sequence of cochain complexes, then $\moc{D^{i},\delta^{i},i\geq 0}$ is called an extension of $\moc{C^{i},\partial^{i},i\geq 0}$ by  $\moc{E^{i},\eta^{i},i\geq 0}$. 
\end{defi}

\begin{defi}
A complex $\moc{C^{i},\partial^{i},i\geq 0}$ is called $0-$\textit{simple} if there exists an integer $n$ such that $\partial^{n}$ is an isomorphism, and $C^{i}$ is trivial for all $i\neq n, n+1$. A complex $\moc{C^{i},\partial^{i},i\geq 0}$ is called $n-$simple if it is an extension of a $(n-1)-$simple complex by another one of type $0-$simple. A complex is called simple if it is $k-$simple for some $k$.
\end{defi}
An easy induction on the height of BG modules yields:
\begin{lem}\label{nsimple}
	If 
	$$0\to C^{0}\xrightarrow{\partial^{0}}C^{1}\xrightarrow{\partial^{1}}\cdots\xrightarrow{\partial^{k}}C^{k+1}\to 0$$
	is an exact sequence of BG modules of finite height, then it is simple.
\end{lem}

In what follows, we characterize minimal resolutions of finite unstable modules.
\begin{prop}
	Let $M$ be a finite unstable module and $\moc{I^{k},\partial^{k},k\geq 0}$ be an injective resolution of $M$. Then this resolution is minimal if and only if every induced morphism 
	$$J(r)\hookrightarrow I^{n}\xrightarrow{\partial^{n}}I^{n+1}\twoheadrightarrow J(s)$$
	is not an isomorphism.
\end{prop}
\begin{proof}
	Suppose that there exist an isomorphism
	$$J(r)\hookrightarrow I^{n}\xrightarrow{\partial^{n}}I^{n+1}\twoheadrightarrow J(r),$$
	then following Lemma \ref{bg2}, the resolution $\moc{I^{k},\partial^{k},k\geq 0}$ is an extension of another resolution of $M$ by an exact $0-$simple sequence and hence cannot be minimal.
	
	Vice versa, let $\moc{J^{k},\delta^{k},k\geq 0}$ be the minimal resolution of $M$. Then there exists a split projection of complexes $\alpha^{i}:I^{i}\to J^{i}$. The kernel of this projection is therefore a finite exact sequence of BG modules. Hence Lemma \ref{nsimple} and the proposition hypothesis imply the triviality of this kernel. Then, $\moc{\alpha^{i},i\geq 0}$ is an isomorphism of complexes.
\end{proof}
\begin{cor}
	Let $M$ be a finite unstable module  and $I:=\moc{I^{k},\partial^{k},k\geq 0}$ be an injective resolution of $M$. Denote by $C$ its maximal simple sub-complex. Then the quotient $I/C$ is the minimal injective resolution of $M$. 
\end{cor}
\begin{defi}\label{iterated}
	Let $I:=\left(I^{k},\partial^{k}\right)$ be a complex of BG modules. We denote by $A(I)$ the quotient $I/C$ of $I$ by its maximal simple sub-complex.
\end{defi}

In the rest of this section, we study an interesting interaction of general exact sequences of BG modules with the \phr procedure. Recall that the Bockstein sequence is defined in Definition \ref{bocksteins}. 
\begin{prop}\label{saturation}
For all integers $n\leq -1$, denote by $\moc{I_{k}^{n,\bullet},\partial^{n,\bullet}}$ the minimal injective resolution of $\Sigma^{k}\mathscr{B}_{n}$. Then, we have an equivalence of exact sequences:
\ali{A\moc{\psr{\Sigma^{k}\mathscr{B}_{\bullet}}{I_{k}^{\bullet,\bullet}}}_{\bullet}\cong \mathscr{B}_{\bullet-k}.}
\end{prop}
\begin{proof}
The case $k=1$ is straight forward as the minimal injective resolution of $\Sigma^{k}\mathscr{B}_{n}$ is given by Mahowald short exact sequences. Suppose that the result is true for all $k\leq m$, we prove that it still holds for $k=m+1$. For all $s\leq -1$, denote by $J^{s,\bullet}$ the minimal injective resolution of $\Sigma A\moc{\psr{\Sigma^{m}\mathscr{B}_{\bullet}}{I_{m}^{\bullet,\bullet}}}_{s}$. We will prove that
\alilabel{suspensioncommute}{A\moc{\psr{\Sigma A\moc{\psr{\Sigma^{m}\mathscr{B}_{\bullet}}{I_{m}^{\bullet,\bullet}}}}{J^{\bullet,\bullet}}}\cong A\moc{\psr{\Sigma^{m+1}\mathscr{B}_{\bullet}}{I_{m+1}^{\bullet,\bullet}}},}
Let $\alpha\leq -1$, then  
\ali{A\moc{\psr{\Sigma A\moc{\psr{\Sigma^{m}\mathscr{B}_{\bullet\geq\alpha}}{I_{m}^{\bullet,\bullet}}}}{J^{\bullet,\bullet}}}\cong A\moc{\psr{\Sigma^{m+1}\mathscr{B}_{\bullet\geq \alpha}}{I_{m+1}^{\bullet,\bullet}}},}
as they are both isomorphic to the minimal injective resolution of $\Sigma^{m+1}\ke{\beta_{\alpha}}.$ Let $\alpha$ tends to $-\infty$, we get the isomorphism \eqref{suspensioncommute}. Therefore, we have:
\ali{A\moc{\psr{\Sigma^{m+1}\mathscr{B}_{\bullet}}{I_{m+1}^{\bullet,\bullet}}}&\cong A\moc{ \psr{\Sigma A\moc{\psr{\Sigma^{m}\mathscr{B}_{\bullet}}{I_{m}^{\bullet,\bullet}}}}{J^{\bullet,\bullet}}},\\
&\cong A\moc{ \psr{\Sigma A\moc{\mathscr{B}_{\bullet-m}}}{J^{\bullet,\bullet}}},\\
&\cong A\moc{ \psr{\Sigma \mathscr{B}_{\bullet-m}}{J^{\bullet,\bullet}}},\\
&\cong \mathscr{B}_{\bullet-m-1}.}
The proposition follows.
\end{proof}

\section{Resolutions of spheres}\label{ressph}
Recall that the Ext-groups $\ext{\u}{s}{\Sigma^{n}\f2}{\Sigma^{t}\f2},$ where $s,t,n\geq 0,$ are of interest because of the unstable Adams spectral sequence converging to the homotopy groups of spheres. One of the most basic way to compute these groups relies on the construction of the minimal injective resolution of $\Sigma^{t}\f2$ for all integers $t\geq 0$. That is what we are going to do in this section. In fact, we will compare the minimal injective resolution of $\Sigma^{t}\f2$ with the Bockstein sequence, and show that they agree in a certain specific area. The keys to this result are Proposition \ref{saturation} and the injective dimension of $\Sigma^{m}J(n)$ for all integers $m,n\geq 0$.

\subsection{Injective dimension}
This section aims at studying the injective dimension of $\Sigma^{m}J(n)$ for all integers $m\geq 0$ and $n\geq 2$. To be more precise, we will show that this dimension is bounded by $\mocv{(m+n)/2}$,  where $[-]$ denotes the integral part of a number. The intrigued reader might wonder why we do not consider the case $n=1$. Here is the reason: this case does not follow the same rule as that of the cases $n\geq 2$. The injective resolution $\mathscr{G}_{m,1}$ of $\Sigma^{m}J(1)$ is of length $m$, and, actually, we will show that $m$ is the injective dimension of $\Sigma^{m}J(1)$.
\begin{lem}\label{injdim2}
For all integers $m\geq 0,n\geq 2$, then the injective dimension of $\Sigma^{m}J(n)$ is bounded by $\mocv{(m+n)/2}$.
\end{lem}
\begin{proof}
We proceed by induction on $m$. The case $m=0$ is trivial. Suppose that we have proved the cases $m<k$, we now verify the case $m=k$. But, it follows from the exact sequence
\ali{0\to \Sigma^{k}J(n)\to\Sigma^{k-1}J(n+1)\to \Sigma^{k-1}J\moc{\frac{n+1}{2}}\to 0}
that 
\ali{\mathrm{dim}_{inj}\Sigma^{k}J(n)=\mathrm{dim}_{inj}\Sigma^{k-1}J(n+1)\leq \mocv{\frac{n+k}{2}}}
if $n$ is even, and that
\ali{\mathrm{dim}_{inj}\Sigma^{k}J(n)&\leq \max\mocn{\mathrm{dim}_{inj}\Sigma^{k-1}J(n+1),\mathrm{dim}_{inj}\Sigma^{k-1}J\moc{\frac{n+1}{2}}}\\
&\leq \max\mocn{\mocv{\frac{n+k}{2}},\mocv{\frac{k-1+\frac{n+1}{2}}{2}}}\leq \mocv{\frac{n+k}{2}}}
if $n$ is odd. This concludes the proof of the lemma.
\end{proof}
\subsection{Resolution of spheres}
In this subsection, we will prove that $\Sigma^{n}\f2$ is of injective dimension $n-1$, and we also describe a large part of its minimal injective resolution.
\begin{thm}
Given an integer $n\geq 0$ and denote by $\moc{I^{k},\partial^{k},k\geq 0}$ the minimal injective resolution of $\Sigma^{n}\f2$. Then, $I^{k}=0$ for all $k>n-1$, and for all $k>\mocv{n/2}$, the resolution coincides with the Bockstein sequence as follows:
\ali{I^{k}\cong \mathscr{B}_{k-n}.}
\end{thm}
\begin{proof}
This is a direct consequence of Proposition \ref{saturation}. Recall that, for all integers $m\leq -1$, we denote by $\moc{I_{k}^{m,\bullet},\partial^{m,\bullet}}$ the minimal injective resolution of $\Sigma^{k}\mathscr{B}_{m}$. Then, we have an equivalence of exact sequences:
\ali{A\moc{\psr{\Sigma^{k}\mathscr{B}_{\bullet}}{I_{k}^{\bullet,\bullet}}}_{\bullet}\cong \mathscr{B}_{\bullet-k}.}
Now, it follows from Lemma \ref{injdim2} that 
\ali{\mathrm{dim}_{inj}\mathscr{B}_{m}\leq \mocv{\frac{k-m}{2}}}
for all $m\leq -1$. Hence, we have
\ali{\mathscr{B}_{s-k} &\cong A\moc{\psr{\Sigma^{k}\mathscr{B}_{\bullet}}{I_{k}^{\bullet,\bullet}}}_{s}\\
&\cong I_{k}^{-1,s}}
for all $s>\mocv{k+1/2}$. Remark that $I^{-1,\bullet}_{k}$ is the minimal injective resolution of $\Sigma^{k}J(1)$, which is isomorphic to $\Sigma^{k+1}\f2$, then replacing $k$ by $n-1$ concludes the proof of the theorem.
\end{proof}
\begin{cor}\label{ressphere}
Let $n,t$ be nonnegative integers and $s> [n/2]$. Then:
     	 	\begin{equation*}
     	 	\ext{\u}{s}{\Sigma^{t}\f2}{\Sigma^{n}\f2}\cong\left\{\begin{array}{ll}
     	 	\f2&\text{ if }t=n-s,\\
     	 	\f2&\text{ if }n-s-1\equiv 0 (4)\text{ and }n-s-1=2t,\\
     	 	0&\text{ otherwise.}
     	 	\end{array}\right.
     	 	\end{equation*}
\end{cor}

\section{Algebraic EHP sequences}\label{ehps}
In this section, we focus on computing extensions groups of unstable modules using the pseudo-hyper resolution method. We then show how to connect these groups by long exact sequences.
\subsection{Extensions groups of unstable modules}
The $J(n)$ form a system of co-generators for the category $\u$. But we can even say more about the classification of injective unstable modules after the work of Lannes and Schwartz. Before going further, recall that there exist other interesting injective unstable modules beside Brown-Gitler modules. These are given by the cohomologies $\mathrm{H}^{*}\left(B\left(\z/2\right)^{\oplus d};\f2\right)$. Lannes and Zarati observe that the tensor products between these modules and Brown-Gitler modules remain injective. Latter, Lannes and Schwartz prove that an indecomposable injective module must be a tensor product of the form $L\otimes J(n)$, where $L$ denotes a certain indecomposable direct summand of $\mathrm{H}^{*}\left(B\left(\z/2\right)^{\oplus d};\f2\right)$. Since the category $\u$ is locally noetherian, we get:
\begin{thm}[\cite{LS89}]
	Each indecomposable injective unstable module is a tensor product of the form $L\otimes J(n)$ between an indecomposable direct summand of a certain  $\mathrm{H}^{*}\left(B\left(\z/2\right)^{\oplus d};\f2\right)$ and a Brown-Gitler module $J(n)$. Every injective unstable module splits into a direct sum between a BG module and a direct sum of modules of the form $L\otimes J(n)$ where $L$ is a direct summand of $\tilde{\mathrm{H}}^{*}\left(B\z/2;\f2\right)^{\otimes d}$ for some $d$. 
\end{thm}
The module $\mathrm{H}^{*}\left(B\z/2;\f2\right)$ has several amazing properties that we will recall now. It is reduced. The functor $-\otimes \mathrm{H}^{*}\left(B\z/2;\f2\right)$ admits a left adjoint in $\u$ which is denoted by $T$. Because $$\mathrm{H}^{*}\left(B\z/2;\f2\right)\cong\tilde{\mathrm{H}}^{*}\left(B\z/2;\f2\right)\oplus \f2$$
then $T\cong \bar{T}\oplus\f2$ where $\bar{T}$ is left adjoint to  $-\otimes \tilde{\mathrm{H}}^{*}\left(B\z/2;\f2\right)$. Moreover, if $M$ is a finite unstable module then $\bar{T}M=0$. Therefore:
\begin{lem}\label{trivialhom}
	Let $M$ be a finite unstable module and $L$ be a direct summand of $\tilde{\mathrm{H}}^{*}\left(B\z/2;\f2\right)^{\otimes d}$ for some $d$. Then
	$$\ho{\u}{M}{L\otimes J(n)}=0$$
	for very $n$.
\end{lem}
\begin{rem}\label{bgpart}
Let $M$ be an unstable module. Denote by $\left(I^{\bullet},\partial^{\bullet}\right)$ its minimal injective resolution. For each $I^{j}$, denote by $B^{j}$ the BG module part and by $R^{j}$ the other. It follows from Lemma \ref{trivialhom} that $\left(B^{\bullet},\partial^{\bullet}\right)$ is a sub-complex of $\left(I^{\bullet},\partial^{\bullet}\right)$. It also follows that if $N$ is a finite unstable module then
$$\ext{\u}{s}{N}{M}\cong \mathrm{H}^{s}\ho{\u}{N}{B^{\bullet}}.$$
\end{rem}
We will now study the extension groups $\ext{\u}{s}{N}{\Sigma^{k}M}$ using the pseudo-hyper resolution method. Let us recall how to use this method to construct an injective resolution of $\Sigma^{k}M$ from  $\left(I^{\bullet},\partial^{\bullet}\right)$. The minimal injective resolution of $\Sigma^{k}J(n)$ is denoted by $B(k,n)$. If $L$ is a direct summand of $\tilde{\mathrm{H}}^{*}\left(B\z/2;\f2\right)^{\otimes d}$ then $L\otimes B(k,n)$ is the minimal injective resolution for $L\otimes \Sigma^{k}J(n)$. Hence there is no BG module part for this resolution. Now using the pseudo-hyper resolution method with respect to these minimal injective resolutions, we get an injective resolution for $\Sigma^{k}M$. Denote this resolution by $M(k)$ and by $MB(k)$ the pseudo-hyper complex obtain from $\Sigma^{k}B^{\bullet}$ using the same method. What we have shown so far is that $MB(k)$ is the BG module part of $M(k)$. Therefore, we have
$$\ext{\u}{s}{N}{\Sigma^{k}M}\cong \mathrm{H}^{s}\ho{\u}{N}{MB(k)}.$$ 
We write $B(k)$ for the quotient of $MB(k)$  by its maximal simple sub-complex. Hence:
$$\ext{\u}{s}{\Sigma^{t}\f2}{\Sigma^{k}M}\cong \ho{\u}{\Sigma^{t}\f2}{B(k)^{s}}.$$ 
We now show how to connect these groups in long exact sequences.
\subsection{Algebraic EHP sequences}
Fix $M$ and $B(t)$ as in the previous sub-section. We index all the direct summands isomorphic to $J(n)$ of $B(t)^{s}$ by the set $A_{(t,n,s)}$, then:
\begin{equation}
\ext{\u}{s}{\Sigma^{n}\f2}{\Sigma^{t}M}\cong \ho{\u}{\Sigma^{n}\f2}{\bigoplus_{A_{(t,n,s)}}J(n)}.\tag{ext}
\end{equation}
From the construction of $\mathscr{G}\moc{B(t)}$ (see Proposition \ref{BGS}), the sum of all the direct summands $J(n)$ of $\mathscr{G}\moc{B(t)}^{s}$ is 
$$B_{(n,t+1,s)}:=C_{(n,t+1,s)}\bigoplus D_{(n,t+1,s)}$$
where 
\begin{align*}
C_{(n,t+1,s)}&:=\bigoplus_{A_{(t,n-1,s)}}J(n),\\ D_{(n,t+1,s)}&:=\bigoplus_{A_{(t,2n-1,s-1)}}J(n).
\end{align*}
Let $X\in \left\{C_{(n,t+1,s)},D_{(n,t+1,s)}\right\}$ and $Y\in \left\{C_{(n,t+1,s+1)},D_{(n,t+1,s+1)}\right\}$, we denote by $\partial_{X,Y}^{s}$ the composition
$$X\hookrightarrow \mathscr{G}\moc{B(t)}^{s}\xrightarrow{\partial^{s}}\mathscr{G}\moc{B(t)}^{s+1}\twoheadrightarrow Y.$$ 
Therefore the maps
$$\partial_{DA}^{s}:=\begin{pmatrix}
\partial_{C_{(n,t+1,s)},C_{(n,t+1,s+1)}}^{s} & \partial_{D_{(n,t+1,s)},C_{(n,t+1,s+1)}}^{s}\\
\partial_{C_{(n,t+1,s)},D_{(n,t+1,s+1)}}^{s} & \partial_{D_{(n,t+1,s)},D_{(n,t+1,s+1)}}^{s}
\end{pmatrix}$$
make $\left\{B_{(n,t+1,s)},\partial_{DA}^{s}\right\}$ a complex. Remark that the BG algorithm shows that
$$\partial_{DA}^{s}:=\begin{pmatrix}
0 & \partial_{D_{(n,t+1,s)},C_{(n,t+1,s+1)}}^{s}\\
0 & 0
\end{pmatrix}.$$
Hence
\begin{align*}
\ext{\u}{s}{\Sigma^{n}\f2}{\Sigma^{t+1}M}&\cong \mathrm{H}_{s}\ho{\u}{\Sigma^{n}\f2}{\mathscr{G}\moc{B(t)}},\\
&\cong \mathrm{H}_{s}\ho{\u}{\Sigma^{n}\f2}{B_{(n,t+1,s)}},\\
&\cong \ke{\left(\partial_{D_{(n,t+1,s)},C_{(n,t+1,s+1)}}^{s}\right)_{*}}\bigoplus \cke{\left(\partial_{D_{(n,t+1,s-1)},C_{(n,t+1,s)}}^{s-1}\right)_{*}}.
\end{align*}
We will now explicit the morphisms $\left(\partial_{D_{(n,t+1,s)},C_{(n,t+1,s+1)}}^{s}\right)_{*}$. The source of this morphism is described as follows:
\begin{align*}
\ho{\u}{\Sigma^{n}\f2}{D_{(n,t+1,s)}}&\cong \ho{\u}{\Sigma^{n}\f2}{\bigoplus_{A_{(t,2n-1,s-1)}}J(n)},\\
&\cong \ho{\u}{\Sigma^{2n-1}\f2}{\bigoplus_{A_{(t,2n-1,s-1)}}J(2n-1)},\\
&\cong \ext{\u}{s-1}{\Sigma^{2n-1}\f2}{\Sigma^{t}M}.
\end{align*}
And the target is:
\begin{align*}
\ho{\u}{\Sigma^{n}\f2}{C_{(n,t+1,s+1)}}&\cong \ho{\u}{\Sigma^{n}\f2}{\bigoplus_{A_{(t,n-1,s+1)}}J(n)},\\
&\cong \ho{\u}{\Sigma^{n-1}\f2}{\bigoplus_{A_{(t,n-1,s+1)}}J(n-1)},\\
&\cong \ext{\u}{s+1}{\Sigma^{n-1}\f2}{\Sigma^{t}M}.
\end{align*}
Therefore 
$$\left(\partial_{D_{(n,t+1,s)},C_{(n,t+1,s+1)}}^{s}\right)_{*}: \ext{\u}{s-1}{\Sigma^{2n-1}\f2}{\Sigma^{t}M}\to \ext{\u}{s+1}{\Sigma^{n-1}\f2}{\Sigma^{t}M}.$$
We denote by $P^{s}$ the morphisms $\left(\partial_{D_{(n,t+1,s)},C_{(n,t+1,s+1)}}^{s}\right)_{*}$ and then we obtain the following exact sequences:
$$0\to\ke{P^{s-1}}\to \ext{\u}{s-2}{\Sigma^{2n+1}\f2}{\Sigma^{t}M}\to \ext{\u}{s}{\Sigma^{n}\f2}{\Sigma^{t}M}\to \cke{P^{s-1}}\to 0,$$
$$0\to \cke{P^{s-1}}\to \ext{\u}{s}{\Sigma^{n+1}\f2}{\Sigma^{t+1}M}\to \ke{P^{s}}\to 0.$$
We can now conclude:
\begin{thm}[Algebraic EHP sequences]\label{EHP}
	There exist a long exact sequence for each $n$:
	\begin{equation}
	\cdots\xrightarrow{H}E_{2}^{s-2,t}(S^{2n+1},M)\xrightarrow{P}E_{2}^{s,t}(S^{n},M)\xrightarrow{E}E_{2}^{s,t+1}(S^{n+1},M)\xrightarrow{H}E_{2}^{s-1,t}(S^{2n+1},M)\xrightarrow{P}\cdots
	\end{equation}
	where $E_{2}^{s,t}(S^{n},M)$ stands for $\ext{\u}{s}{\Sigma^{n}\f2}{\Sigma^{t}M}$. 
\end{thm}

\subsection{A special case}
This section studies the special case of EHP sequence for $n=2^{k}$. In particular we show that in this case the morphisms $P$ are trivial. 
\begin{lem}
	If $n=2^{k}$ then the morphism $\left(\partial_{D_{(n,t+1,s)},C_{(n,t+1,s+1)}}^{s}\right)_{*}$ is trivial. 	
\end{lem}
\begin{proof}
	Recall that 
	$$\partial_{D_{(2^{k},t+1,s)},C_{(2^{k},t+1,s+1)}}^{s}:\bigoplus_{A_{(t,2^{k+1}-1,s-1)}}J(2^{k})_{\omega}\to \bigoplus_{A_{(t,2^{k}-1,s+1)}}J(2^{k})_{\alpha}.$$
	Suppose that this morphism is not trivial then there must exist an identity coefficient in the matrix form. The first item of the BG algorithm assures that in this case, the Steenrod square $Sq^{2^{k}}$ can be factorized as a product of other Steenrod operations. This cannot be true since $Sq^{2^{k}}$ is indecomposable. 	
	Therefore $\partial_{D_{(2^{k},t+1,s)},C_{(2^{k},t+1,s+1)}}^{s}$ is trivial.
\end{proof}
We get the following consequence:
\begin{thm}[James splitting]\label{jamessplit}
	There are short exact sequences
	$$0\to E^{s,t}_{2}\left(S^{2^{k}-1},M\right)\to E^{s,t+1}_{2}\left(S^{2^{k}},M\right)\to E^{s-1,t}_{2}\left(S^{2^{k+1}-1},M\right)\to 0$$
	for every $k$.
\end{thm}
\section{Bousfield's proof on the existence of algebraic EHP sequences}\label{BS}
\begin{flushright}
	\begin{minipage}[r]{29em}\fontsize{9pt}{9pt}\selectfont
		\textit{An interesting fact about the algebraic EHP sequence: it can be derived in a completely abstract way. That is, it can be derived without the construction of special projective or injective resolutions and without any computation whatsoever. Bousfield explained to me how to do this, about 45 years ago. Here is the key idea. One has a "loop functor" on the category of unstable Steenrod modules. It is left adjoint to the suspension. This functor is right exact, and has non-trivial left-derived functors. The key is to notice that these left-derived functors are zero, in homological degrees greater than one. The existence of the long-exact EHP sequence follows immediately.}\vspace{1em}
		
		\hfill{}{\textsc{William M. Singer}, \textit{Private communication }\cite{Wil16}}\vspace{1em}
	\end{minipage}
\end{flushright}
Let $M$ and $N$ be two unstable modules.  In this section, we prove that there exists a long exact sequence of Ext-groups
$$\ext{\u}{s-2}{\Omega_{1}M}{N}\xrightarrow{}\ext{\u}{s}{\Omega M}{N}\xrightarrow{}\ext{\u}{s}{M}{\Sigma^{}N}\xrightarrow{}\ext{\u}{s-1}{\Omega_{1}M}{N}$$
for all integers $s\geq 2$, where $\Omega$ is the left adjoint of the suspension functor $\Sigma$ and $\Omega_{k}$ is the $k-$th left-derived functor of $\Omega$. 
\subsection*{The loop functor of unstable modules}
Recall that the Frobenius twist of unstable modules exists under the form of the double functor $\Phi$, defined in Section \ref{section2}. The operator $\mathrm{Sq}_{0}$, which associates an element $x$ of degree $n$ of a certain unstable module $M$ with $Sq^{n}x$, defines an $\a2-$linear morphism $\lambda_{M}:\Phi M\to M$. It is classical that both kernel and cokernel of this morphism are a suspension of an unstable module, which can be described explicitly according to the following proposition.
\begin{prop}
	Let $M$ be an unstable module. Then, the cokernel of $\lambda_{M}$ is a suspension, and we denote by $\Omega M$ the unstable module such that $\Sigma\Omega M\cong \cke{\lambda_{M}}$. Then, $\Omega$ defines an endofunctor of unstable modules which is also a left adjoint of $\Sigma$. Let $\Omega_{k}M$ be the $k-$left-derived functor of $\Omega$ for all integers $k\geq 0$, then we have
	\ali{\Sigma\Omega_{k}M\cong\left\{ \begin{array}{cl}
	\ke{\lambda_{M}}&\text{ if }k=1,\\
	0&\text{ if }k>1. 
	\end{array}\right.} 
\end{prop}
\begin{proof}
An unstable module $N$ is isomorphic to $\Sigma Q$ for some $Q\in\u$ if and only if $\mathrm{Sq}_{0}$ acts trivially on $N$. Then, it is evident that the cokernel of $\lambda_{M}$ is a suspension of some unstable module that we denote by $\Omega M$. This is a well-defined endofunctor of $\u$ as $\lambda_{M}$ is natural on $M$. And since both $\Phi$ and $Id$ are exact, it follows that $\Omega$ is right exact. In order to show that $\Omega$ is left adjoint to $\Sigma$, it suffices to show that 
\ali{\ho{\u}{\Omega P}{N}\cong \ho{\u}{P}{\Sigma N}}
for all unstable modules $P,N$ such that $P$ is projective. But this is a direct consequence of the following well-known short exact sequence (see, e.g., \cite{Sch94}).
\ali{
0\to \Phi F(n)\xrightarrow{\lambda(F(n))}F(n)\xrightarrow{\sigma(F(n))}\Sigma F(n-1)\to 0.}
where $\sigma\moc{F(n)}$ is defined by $\imath_{n}\mapsto\Sigma\imath_{n-1}$.

Now, to study the left-derived functors of $\Omega$, let $\left(P_{i},\partial_{i}:P_{i+1}\to P_{i}\right)_{i\geq 0}$, abbreviated as $P_{\bullet}$, be a projective resolution of $M$. Because, on the one hand, $\Sigma \Omega P_{\bullet}$ fits in the short exact sequence of complexes
$$0\to \Phi P_{\bullet}\xrightarrow{\lambda_{P_{\bullet}}} P_{\bullet}\xrightarrow{\sigma_{P_{\bullet}}} \Sigma \Omega P_{\bullet}\to 0,$$
and, on the other hand, $\Phi$ is exact, then the homology groups of $\Sigma \Omega P_{\bullet}$ are trivial, in degrees greater than one, and there is a short exact sequence,
\begin{equation}\label{unit}
0\to \mathrm{H}_{1}\left(\Sigma \Omega P_{\bullet}\right)\to \Phi M\xrightarrow{\lambda_{M}} M\xrightarrow{\sigma_{M}} \Sigma\Omega M\to 0,
\end{equation}
connecting the homology groups of degrees $0$ and $1$. As $\Sigma$ is exact, we have:
$$\ke{\lambda_{M}}\cong \Sigma\Omega_{1}M.$$	
The proposition follows.
\end{proof}
\begin{rem}\label{trivial}
	\begin{itemize}
		\item For all unstable module $M$, the morphism $\lambda_{\Sigma M}$ is trivial. Therefore, $$\Omega\Sigma M\cong M,\quad\forall M\in \u.$$
		\item The loop of $\sigma_{M}$ is the identity of $\Omega M$. As the loop functor $\Omega$ is right exact, then $\Omega \lambda_{M}$ is trivial. 
		\item Fix $\left\{P_{i},\partial_{i}:P_{i+1}\to P_{i},i\geq 0\right\}$ a projective resolution of $M$. Denote by $C$ the co-kernel $\cke{\Omega\partial_{1}}$. Then $\left\{\Omega P_{i},\Omega\partial_{i}:\Omega P_{i+1}\to \Omega P_{i},i\geq 1\right\}$  is a projective resolution of $C$. Moreover, $C$ fits in the short exact sequence:
		$$0\to \Omega_{1}M\to C\to \frac{\Omega P_{1}}{\ke{\Omega\partial_{0}}}\to 0.$$
		Because $\Omega$ is right exact, $\Omega P_{1}/\ke{\Omega\partial_{0}}$ is isomorphic to the kernel of the morphism $\Omega P_{0}\twoheadrightarrow \Omega M.$ 
		\item If $\Omega_{1}M$ is trivial, then $\Omega P_{\bullet}$ is a projective resolution of $\Omega M$.
	\end{itemize}	
\end{rem}
\subsection*{Bousfield's sequences}
Let $N$ be an unstable module. Let $M$ be an unstable module and $\left\{P_{i},\partial_{i}:P_{i+1}\to P_{i},i\geq 0\right\}$, abbreviated as $P_{\bullet}$, be a projective resolution of $M$. Since $\Omega P_{0}$ is projective, the long exact sequence of Ext-groups associated with the short exact sequence
$$0\to  \frac{\Omega P_{1}}{\ke{\Omega\partial_{0}}}\to \Omega P_{0}\to \Omega M\to 0$$ splits into an exact sequence
$$0\to \ho{\u}{\Omega M}{N}\to \ho{\u}{\Omega P_{0}}{N}\to \ho{\u}{ \frac{\Omega P_{1}}{\ke{\Omega\partial_{0}}}}{N}\to \ext{\u}{1}{\Omega M}{N}\to 0,$$
and isomorphisms
$$\ext{\u}{s}{ \frac{\Omega P_{1}}{\ke{\Omega\partial_{0}}}}{N}\xrightarrow{\sim} \ext{\u}{s+1}{\Omega M}{N},$$
for all $s\geq 1$. Now, because  $\left\{\Omega P_{i},\Omega\partial_{i}:\Omega P_{i+1}\to\Omega P_{i},i\geq 1\right\}$ is a projective resolution of $C$ (see Remark \ref{trivial}), then for every $s\geq 1$ we have:
\begin{align*}
\ext{\u}{s}{C}{N}&\cong \mathrm{H}^{s+1}\left(\ho{\u}{\Omega P_{\bullet}}{N},\left(\Omega\partial_{\bullet}\right)^{*}\right)\\
&\cong \ext{\u}{s+1}{M}{\Sigma N}.
\end{align*}
Therefore, the long exact sequence of Ext-groups associated with the short exact sequence 
$$0\to \Omega_{1}M\to C\to \frac{\Omega P_{1}}{\ke{\Omega\partial_{0}}}\to 0$$
is the general algebraic long-exact EHP sequence.
\begin{thm}[Bousfield's proof]\label{Bousfield}
	For all unstable modules $M$ and $N$, there exists a long exact sequence
	$$\cdots\xrightarrow{}\ext{\u}{s-2}{\Omega_{1}M}{N}\xrightarrow{}\ext{\u}{s}{\Omega M}{N}\xrightarrow{}\ext{\u}{s}{M}{\Sigma^{}N}\xrightarrow{}\ext{\u}{s-1}{\Omega_{1}M}{N}\xrightarrow{}\cdots$$
	of Ext-groups, where $s\geq 2$.
\end{thm}	
Let $M$ be $\Sigma^{n}\f2$ and $N$ be $\Sigma^{t}\f2$. If $n\geq 1$, then the morphism $\lambda_{M}:\Phi M\to M$ is trivial. Therefore,
\ali{
\Omega M&\cong \Sigma^{n-1}\f2,\\
\Omega_{1} M&\cong\Sigma^{2n-1}\f2.
}
A reformulation of Bousfield's long exact sequence, in this case, yields the algebraic EHP sequence for $S^{n}$.
\begin{thm}
	For every positive integer $n$, there exists a long exact sequence
	\begin{equation*}
			\cdots\xrightarrow{H}E_{2}^{s-2,t}(S^{2n+1})\xrightarrow{P}E_{2}^{s,t}(S^{n})\xrightarrow{E}E_{2}^{s,t+1}(S^{n+1})\xrightarrow{H}E_{2}^{s-1,t}(S^{2n+1})\xrightarrow{P}\cdots,
	\end{equation*}
	where $E_{2}^{s,t}(S^{n}):=\ext{\u}{s}{\Sigma^{n}\f2}{\Sigma^{t}\f2}$.
\end{thm}
\subsection*{Application}
In this subsection, we use the loop functor $\Omega$ to study a special case of the algebraic EHP sequence. 

If $\left\{C_{i},\partial_{i}:C_{i+1}\to C_{i},i\geq 0\right\}$ is a complex, denote by $C_{\bullet}[1]$ the complex:
\begin{align*}
C_{\bullet}[1]_{i}&=\left\{ \begin{array}{cl}
C_{i-1}&\textup{ if }i\geq 1,\\
0&\textup{ if }i=0,
\end{array} \right.\\
\partial[1]_{i}&=\left\{ \begin{array}{cl}
\partial_{i-1}&\textup{ if }i\geq 1,\\
0&\textup{ if }i=0.
\end{array} \right. 
\end{align*}
Let $M$ be an unstable module such that $\Omega_{1}M$ is trivial. Fix $\left\{P_{\bullet},\partial_{i}:P_{i+1}\to P_{i},i\geq 0\right\}$, abbreviated as $P_{\bullet}$, a projective resolution of $M$, and fix $\left\{Q_{\bullet},\delta_{i}:Q_{i+1}\to Q_{i},i\geq 0\right\}$, abbreviated as $Q_{\bullet}$, a projective resolution of $\Phi M$. The natural transformation $\lambda:\Phi\to Id$ gives rise to a morphism of complexes: $\lambda_{P_{\bullet}}:\Phi P_{\bullet}\to P_{\bullet}$. On the other hand, the identity of $\Phi M$ yields a morphism of complexes: $\omega:Q_{\bullet}\to \Phi P_{\bullet}$. Therefore, the composition map $\omega\circ\lambda_{P_{\bullet}}$ makes the following diagram commute.
$$\xymatrixcolsep{4em}\xymatrixrowsep{1.5em}
\xymatrix{Q_{\bullet}\ar[r]^{\omega\circ\lambda_{P_{\bullet}}}\ar[d]&P_{\bullet}\ar[d]\\
	\Phi M\ar[r]_{\lambda_{M}}&M
}
$$
Now, we can consider $\omega\circ\lambda_{P_{\bullet}}:Q_{\bullet}\to P_{\bullet}$ as a double complex with two non-trivial columns $Q_{\bullet}$ and $P_{\bullet}$. Denote by $T_{\bullet}$ the total complex of this double complex. As these columns, $Q_{\bullet}$ and $P_{\bullet}$, are acyclic, the homology groups of $T_{\bullet}$ are computed as follows.
$$\mathrm{H}_{i}\left(T_{\bullet}\right)\cong \left\{ \begin{array}{cl}
\cke{\lambda_{M}}&\textup{ if }i=0,\\
\ke{\lambda_{M}}&\textup{ if }i=1,\\
0&\textup{ otherwise}.
\end{array}\right.$$
Since $\Omega_{1}M$ is trivial, $T_{\bullet}$ is a projective resolution of $\Sigma\Omega M$. We now compute $\Omega T_{\bullet}$. It follows from Remark \ref{trivial} that the morphism $\Omega \left(\omega\circ \lambda_{P_{\bullet}}\right)$ is trivial. We then have:
$$\Omega T_{\bullet}\cong \Omega P_{\bullet}\bigoplus \Omega Q_{\bullet}[1].$$
We also deduce from Remark \ref{trivial} that $ \Omega P_{\bullet}$ is a projective resolution of $\Omega M$, and $\Omega Q_{\bullet}$ is a projective resolution of $\Omega \Phi M$. 
\begin{lem}\label{bousfieldlem}
	Let $M$ be an unstable module such that $\Omega_{1}M$ is trivial. For all unstable module $N$, we have an isomorphism of Ext-groups
	$$\ext{\u}{s}{\Sigma\Omega M}{\Sigma N}\cong \ext{\u}{s}{\Omega M}{\Sigma N}\bigoplus \ext{\u}{s-1}{\Omega \Phi M}{N},$$
	for all $s\geq 0$. (Here, by convention, the Ext-groups of degree $-1$ are trivial.) 
\end{lem}
\begin{proof}
	The Ext-groups $\ext{\u}{*}{\Sigma\Omega M}{\Sigma N}$ can be computed as follows.
	\begin{align*}
	\ext{\u}{s}{\Sigma\Omega M}{\Sigma N}&\cong \mathrm{H}^{s}\left(\ho{\u}{T_{\bullet}}{\Sigma N}\right)\\
	&\cong \mathrm{H}^{s}\left(\ho{\u}{\Omega T_{\bullet}}{N}\right)\\
	&\cong \mathrm{H}^{s}\left(\ho{\u}{\Omega P_{\bullet}\bigoplus \Omega Q_{\bullet}[1]}{N} \right)\\
	&\cong  \ext{\u}{s}{\Omega M}{N}\bigoplus \ext{\u}{s-1}{\Omega \Phi M}{N}.
	\end{align*}
	We can then conclude the lemma.
\end{proof}
Remark that, following the exact sequence 
$$\Phi^{2}M\to \Phi M\to \Phi\left(\Sigma\Omega M\right)\to 0,$$
we have 
$$\Omega \Phi M\cong \Sigma^{}\Phi\Omega M$$
for all unstable modules $M$. Then, applying Lemma \ref{bousfieldlem} to $M=\Phi^{n}F(1)$, we recover Theorem \ref{jamessplit}.
\section{Infinite complex projective space}\label{icps}
This section aim at studying a particular relation between the infinite complex projective and spheres, which is a direct consequence of the following result in view of pseudo-hyperresolution method.
\begin{lem}
Denote by $\h$ the cohomology $\h^{*}\moc{B\z/2;\f2}$ of the classifying space of $\z/2$. Then, we have
\ali{\Phi \h&\cong \h^{*}\moc{\mathbb{C}P^{\infty};\f2},\\
\Omega \h&\cong \Phi \h.}
\end{lem}
\begin{proof}
The lemma follows from the fact that $\h$ is isomorphic to the polynomial algebra $\f2\mocv{u}$ on one variable of degree $1$ whereas $\h^{*}\moc{\mathbb{C}P^{\infty};\f2}$ is the polynomial algebra $\f2\mocv{t}$ on one variable of degree $2$.
\end{proof}
\begin{rem}\label{rem112}
Denote by $\bar{\h}$ the reduced cohomology $\tilde{\h}^{*}\moc{B\z/2;\f2}$. Then, we have two short exact sequences
\ali{\xymatrixrowsep{1ex}\xymatrix{
0\ar[r]& \Phi\bar{\h}\ar[r]&\bar{\h}\ar[r]&\Sigma\Phi \h\ar[r] & 0,\\
0\ar[r]& \Phi \h\ar[r]&\h\ar[r]&\Sigma\Phi \h\ar[r] & 0.
}}
As a result, we have the following long exact sequences
\ali{\xymatrixrowsep{1ex}\xymatrix{
0\ar[r]& \Phi\bar{\h}\ar[r]&\bar{\h}\ar[r]&\Sigma \h \ar[r] & \Sigma^{2} \h \ar[r] & \cdots \ar[r]& \Sigma^{d} \h\ar[r] & \cdots,\\
0\ar[r]& \Phi \h\ar[r]&\h\ar[r]&\Sigma \h \ar[r] & \Sigma^{2} \h \ar[r] & \cdots \ar[r]& \Sigma^{d} \h\ar[r] & \cdots.
}}
\end{rem}
\begin{thm}\label{compsphfinal}
For all integers $s\geq 0$, we have
\alilabel{comproj}{
\ext{\u}{s}{\Sigma^{n}\f2}{\Sigma^{k}\Phi \h}\cong\bigoplus_{m+q=s} \ext{\u}{q}{\Sigma^{n}\f2}{\Sigma^{k+m}\f2}.
}
\end{thm}
\begin{proof}
Following Remark \ref{rem112}, there is a long exact sequence
\ali{\xymatrix{
0\ar[r]& \Sigma^{k}\Phi \h\ar[r]&\Sigma^{k} \h\ar[r]&\Sigma^{k+1} \h \ar[r] & \Sigma^{k+2} \h \ar[r] & \cdots \ar[r]& \Sigma^{k+d} \h\ar[r] & \cdots
}}
for all integers $k\geq 0$. In accordance with Section \ref{ehps}, we denote by $B\moc{s-1,1}$ the minimal injective resolution of $\Sigma^{s}\f2$ for all integers $s\geq 0$. As $\h\otimes J(n)$ is an injective unstable module for all integers $n\geq 0$, then $\h\otimes B\moc{s-1,1}$ is the minimal injective resolution of $\Sigma^{s}\h$ for all integers $s\geq 0$. Let $B\moc{s-1,1}^{t}$ denote the $t-$th term of the resolution $B\moc{s-1,1}$, and let $\partial^{l,t}$ be the differential 
\ali{\h\otimes B\moc{k-1+l,1}^{t}\to
	\h\otimes B\moc{k-1+l,1}^{t+1}.}
Then, in view of Proposition \ref{key12}, there exist morphisms
	$$\partial^{l,t}_{i}:\h\otimes B\moc{k-1+l,1}^{t}\to
	\h\otimes B\moc{k-1+l+i+1,1}^{t-i}$$ 
	such that the morphisms
	\alilabel{partial}{\partial^{l}=\begin{pmatrix}
	\partial^{0,l} & 0 & \cdots & 0 \\
	\partial^{0,l}_{0} & \partial^{1,l-1} & \cdots & 0 \\
	\vdots  & \vdots  & \ddots & \vdots  \\
	\partial^{0,l}_{l-1} & \partial^{1,l-1}_{l-2} & \cdots & \partial^{l,0}\\
	\partial^{0,l}_{l} & \partial^{1,l-2}_{l-1} & \cdots & \partial^{l,0}_{0}
	\end{pmatrix}}
	make the sequence
	\begin{equation}\label{res113}
	\h\otimes B\moc{k-1,1}^{0}\xrightarrow{\partial^{0}} \h\otimes B\moc{k-1,1}^{1}\bigoplus
	\h\otimes B\moc{k,1}^{0}\xrightarrow{\partial^{1}}\cdots\xrightarrow{\partial^{l}}\dsum{m+n=l+1}{}{\h\otimes B\moc{k-1+m,1}^{n}}
	\xrightarrow{\partial^{l+1}}\cdots
	\end{equation}
	an injective resolution of $\Sigma^{k}\Phi \h$. It follows from Remark \ref{bgpart} that 
	\ali{\ext{\u}{s}{\Sigma^{n}\f2}{\Sigma^{k}\Phi \h}\cong \ext{\u}{s}{\Sigma^{n}\f2}{B^{\bullet}},}
	where $B^{\bullet}$ is the BG part of the resolution \eqref{res113}. It is evident that 
	\ali{B^{s}\cong \dsum{m+n=s}{}{ B\moc{k-1+m,1}^{n}}.}
	By abuse of notation, we adopt the notation \eqref{partial} for the induced differentials of $B^{\bullet}$. Because a morphism between two BG modules cannot factorize via a direct sum of tensor products of the form $\bar{\h}\otimes J(n)$, then the induced differential $\partial^{l}$ of $B^{\bullet}$ has the following simple form
	\ali{\partial^{l}=\begin{pmatrix}
		\partial^{0,l} & 0 & \cdots & 0 \\
		0 & \partial^{1,l-1} & \cdots & 0 \\
		\vdots  & \vdots  & \ddots & \vdots  \\
		0 & 0 & \cdots & \partial^{l,0}\\
		0 & 0 & \cdots & 0
		\end{pmatrix}}
   Therefore, $B^{\bullet}$ is the direct sum of the resolution $B\moc{s,1}$ with some appropriate shifting. This conclude the isomorphism \eqref{comproj}.	
\end{proof}

\section*{Acknowledgements}
This research was supported through the programme "\textbf{Oberwolfach Leibniz Fellows}" by the Mathematisches Forschungsinstitut Oberwolfach (MFO) in 2016. It is a great pleasure for the author to express his sincere thanks to all members of the MFO for their hospitality during his three-month visit in Oberwolfach, where the present work was accomplished. The author enjoyed various discussions with Professor Singer about the algebraic EHP sequence and would like to thank him for pointing out Bousfield's proof and for providing references concerning the Lambda algebra. The acknowledgements also go to Professor Bousfield for his generous permission for exhibiting his approach to the construction of the algebraic EHP sequence. Part of this work was conducted during the author's PhD thesis \cite{Cuo14b} at University Paris 13, and he is grateful to his advisor, Professor Schwartz, for his guidance and his constant supports. The author also takes advantage of this occasion to thank Professor Touz\'{e} for pointing out to him the similarity between the approach of the present paper to functor homology and the original ones in \cite{FS97} and \cite{Tot97}, and for showing him the existence of the mapping telescope technique \cite{BN93}. 

\bibliographystyle{plain}
\bibliography{MathScinet.bib}
\end{document}

%% file: exam.tex
\begin{tikzpicture}[line cap=round,line join=round,>=triangle 45,x=4cm,y=0.8cm]
\clip(0.5,2.15) rectangle (3.5,6.8);

\node[circle,scale=0.3,fill=black,label={[label distance=0.1cm]180:$J(7)$}] (1a7) at (1, 5) [draw] {};
\node[circle,scale=0.3,fill=black,label={[label distance=0.1cm]180:$J(6)$}] (1a6) at (1, 4) [draw] {};
\node[circle,scale=0.3,fill=black,label={[label distance=0.1cm]0:$J(6)$}] (3a6) at (3, 6) [draw] {};
\node[circle,scale=0.3,fill=black,label={[label distance=0.1cm]0:$J(4)$}] (3a4) at (3, 4) [draw] {};
\node[circle,scale=0.3,fill=black,label={[label distance=0.1cm]0:$J(3)$}] (3a3) at (3, 3) [draw] {};

\draw [->,>=stealth] (1a7) -- node[pos=0.5,fill=white,rounded corners, sloped] {$\scriptstyle\bullet Sq^{1}$} (3a6);
\draw [->,>=stealth] (1a7) -- node[pos=0.5,fill=white,rounded corners, sloped] {$\scriptstyle\bullet Sq^{2,1}$} (3a4);
\draw [->,>=stealth] (1a6) -- node[pos=0.5,fill=white,rounded corners, sloped] {$\scriptstyle\bullet Sq^{2}$} (3a4);
\draw [->,>=stealth] (1a6) -- node[pos=0.5,fill=white,rounded corners, sloped] {$\scriptstyle\bullet Sq^{3}$} (3a3);
\draw \boundellipse{0.9,4.5}{0.2}{1.3};
\draw \boundellipse{3.1,4.5}{0.25}{2.3};
\end{tikzpicture}

%% file: sigma2.tex
\begin{tikzpicture}[line cap=round,line join=round,>=triangle 45,x=3cm,y=1.2cm]
\clip(0.5,-0.5) rectangle (4.5,2.5);

\node[draw=none] (1a2) at (1, 2) {$\Sigma^{2}J(1)$};
\node[draw=none] (2a2) at (2, 2) {$\Sigma^{}J(2)$};
\node[draw=none] (3a2) at (3, 2) {$\Sigma^{}J(1)$};
\node[draw=none] (4a2) at (4, 2) {$0$};

\node[draw=none] (2a1) at (2, 1) {$J(3)$};
\node[draw=none] (3a1) at (3, 1) {$J(2)$};

\node[draw=none] (2a0) at (2, 0) {$0$};
\node[draw=none] (3a0) at (3, 0) {$J(1)$};

\draw [->,>=stealth] (1a2) -- (2a2);
\draw [->,>=stealth] (2a2) -- node[above] {$\scriptstyle\Sigma\moc{\bullet Sq^{1}}$} (3a2);
\draw [->,>=stealth] (3a2) -- (4a2);

\draw [->,>=stealth] (2a1) -- node[above] {$\scriptstyle\bullet Sq^{1}$} (3a1);
\draw [->,>=stealth] (2a0) -- (3a0);

\draw [->,>=stealth] (2a2) -- (2a1);
\draw [->,>=stealth] (3a2) -- (3a1);
\draw [->,>=stealth] (2a1) -- (2a0);
\draw [->,>=stealth] (3a1) -- (3a0);

\draw(1.5,-0.3) rectangle (3.5,1.5);
\end{tikzpicture}